\documentclass[aap]{imsart}

\usepackage{latexsym}
\usepackage{amstext}
\usepackage{tikz}
\usepackage{epsfig,amsmath,amssymb,amscd,amsfonts,bm,subfigure,color,graphicx,color}

\startlocaldefs
 
\newcommand{\ds}{\displaystyle}
\newcommand{\beq}{\begin{eqnarray}}
\newcommand{\eeq}{\end{eqnarray}}
\newcommand{\beqq}{\begin{eqnarray*}}
\newcommand{\eeqq}{\end{eqnarray*}}

\newcommand{\mn}{(m_1,\dots, m_N)}
\newcommand{\ANK}{P_{N,K}}
\newcommand{\ANKM}{P_{N,K,M}}
 \newtheorem{proposition}{Proposition}
\newtheorem{lemma}{Lemma}
\newtheorem{thm}{Theorem}[section]

 \newtheorem{remark}{Remark}
\newcommand{\Proof}{\underline{\bf Proof:}}
\newcommand{\Endproof}{~$\Box$~}

\font\bb=msbm10 at 12pt

\def\rN{\hbox{\bb N}}

\endlocaldefs

\begin{document}

\begin{frontmatter}

\title{Stochastic coagulation-fragmentation processes with a finite number of particles and applications}
\runtitle{Stochastic coagulation-fragmentation }


\author{\fnms{Nathanael} \snm{Hoze} \thanksref{m1} \ead[label=e1]{nathanael.hoze@env.ethz.ch}}
\address{\printead{e1}}
 \and
 \author{\fnms{David} \snm{Holcman} \thanksref{m3,m4,t1}\ead[label=e2]{david.holcman@ens.fr}}
\address{\printead{e2}}

\affiliation{\thanksmark{m1} Institute of Integrative Biology, ETH Zurich, Zurich, Switzerland.
  \thanksmark{m3} Newton Institute and Department of Applied Mathematics and Theoretical Physics (DAMPT), University of Cambridge, CB30DS United Kingdom and
  \thanksmark{m4} Ecole Normale Sup\'erieure, 46 rue d'Ulm 75005 Paris, France}
\thankstext{t1}{David Holcman is supported by a Marie Curie Award and a Simons fellowship. He thanks the Institute of Mathematics, Oxford University, for hospitality.}

\runauthor{N. Hoze and D. Holcman}

\begin{abstract}
{Coagulation-fragmentation processes describe the stochastic association and dissociation of particles in clusters. 
Cluster dynamics with cluster-cluster interactions for a finite number of particles has recently attracted attention especially in stochastic analysis and statistical physics of cellular biology, as novel experimental data is now available, but their interpretation remains challenging.} We derive here probability distribution functions for clusters that can either aggregate upon binding to form clusters of arbitrary sizes or a single cluster can dissociate into two sub-clusters. Using combinatorics properties and Markov chain representation, we compute steady-state distributions and moments for the number of particles per cluster in the case where the coagulation and fragmentation rates follow a detailed balance condition. We obtain explicit and asymptotic formulas for the cluster size and the number of clusters in terms of hypergeometric functions. To further characterize clustering, we introduce and discuss two mean times: one is time two particles spend together before they separate and other is the time they spend separated before they meet again for the first time. Finally we discuss applications of the present stochastic coagulation-fragmentation framework in cell biology.
\end{abstract}

\begin{keyword}[class=MSC]
\kwd[Primary ]{60J20}
\kwd{Markov chain}
\kwd{Coagulation-fragmentation processes}
\kwd{Stochastic processes}
 \kwd[; secondary ]{05A17.}
\end{Integer partition}
\end{frontmatter}

\section*{Introduction}
Clustering appears in various areas of science such as astrophysics, where masses can form aggregate under gravitation, biochemistry where molecules have to meet to react, colloids that aggregate in solution or ecology where prey-predator have to meet to stabilize populations. A century ago, Von Smoluchowski \cite{Smoluchowski1916} described an irreversible aggregation of many particles in clusters.  Later on, a set of coagulation-fragmentation equations was proposed by Becker-D\"oring when clusters can lose or gain only one particle at a time \cite{Becker,Chandrasekar,redner,Wattis}. Nowadays, continuous limit analysis \cite{Amann,Norris}, determinist, stochastic, asymptotical and numerical methods are developed to study clustering, estimate the number of clusters \cite{Aldous,benavraham1988,Thomson1989,Rotstein2015,Collet2004} and their sizes.

{The coagulation models mentioned above have been extended for cluster growth with a finite number of particles \cite{Lushnikov1978, Marcus1968}, known as Marcus-{ Lushnikov} process, and more recently for discrete coagulation-fragmentation models with an infinite number of particles   \cite{Ball1990}.
} Many statistical and probabilistic studies have been made on the Marcus-Lushnikov process \cite{Jacquot2009}, however, much less is known about the statistical properties for the coagulation-fragmentation of a finite number of particles  \cite{Gueron1998}. Of particular interest in stochastic biology are models of coagulation-fragmentation where the cluster size cannot exceed a given threshold \cite{Hoze2012,Yvinec2012}. These models are used in genetics for describing the organization in clusters of the chromosome ends \cite{Hoze2012} or to model viral capsid assembly in cells \cite{Zlotnick2005,Hoze2014,Hoze2015}. This article presents several exact and asymptotic results and it aims to attract attention toward developments of applied probability with direct applications to interpret data.

We recall that the Smoluchowski equations for coagulation-fragmentation consist of an infinite system of differential equations for the number $n_j(t)$ of clusters of size $j$ at time $t$ in a population of infinite size \cite{Smoluchowski1916}. The index $j$ can take values between 1 and $\infty$ and
{\beq
\frac{dn_j(t)}{dt} &=&    \frac{1}{2}\sum_{k=0}^{j-1}A(k,j-k)n_k(t)n_{j-k}(t)   -n_j(t) \sum_{k=1}^{\infty}A(j,k) n_k(t) \nonumber \\
&-&\frac{1}{2}n_j(t) \sum_{k=1}^{j-1}B(k,j-k)   + \sum_{k=1}^{\infty}B(j,k)n_{k+j}(t),
\label{deterministicSmoluchowski}
\eeq}
where the first line in the left-hand side corresponds to the coagulation and the second accounts for the fragmentation. The coagulation kernel $A(i,j)$ is the rate at which two clusters of size $i$ and $j$ coalesce to form a cluster of size $i+j$, while the fragmentation kernel $B(i,j)$ is the rate at which a cluster of size $i+j$ dissociates into a cluster of size $i$ and a cluster of size $j$. { When cluster sizes cannot exceed a threshold $M$,  the system of equations \eqref{deterministicSmoluchowski} is truncated and the coagulation kernel is $A(i,j)=0$ if $i+j\geq M$ \cite{DaCosta1998,Carr1994}. This system of equations is a mean-field deterministic model of the coagulation-fragmentation process that  describes a discrete number of particles aggregating and dissociating, but it does not allow a complete analysis of the cluster distribution in the small size limit \cite{Gueron1998} }.

{ The goal of this paper is to study a stochastic system of coagulation-fragmentation with a limited number of particles.} We present exact solutions and expressions for the distribution, statistical moments of clusters and number of particles per cluster, using the following generic rules of coagulation-fragmentation: a cluster of size $i+j$ can give rise to two clusters of size $i$ and $j$ at a rate $F(i,j)$, and two clusters of size $i$ and $j$ form a new cluster of size $i+j$ at a rate $C(i,j)$.

We focus our analysis on coagulation-fragmentation processes (CFP) that verify the detailed balance condition \cite{Durrett1999}, for which there exists a function $a(i)=a_i $ such that $\forall i,j \in \rN$
\beq
\frac{C(i,j)}{F(i,j)} = \frac{a(i+j)}{a(i)a(j)}.
\eeq
{The detailed balance condition ensures that, at equilibrium, each elementary coagulation or fragmentation process is  equilibrated by its reciprocal process.}
Using exact { formulas} for probabilities for these configurations when the total number of clusters is fixed, we compute the probability distribution function of the number of clusters. We also compute the probability distribution that the number of cluster of size $i$ is $m_i$ so that the distribution of sizes of the ensemble of clusters is $(m_1,\dots,m_n)$. When there are exactly $N$ particles and the total number of clusters is fixed to $K$, we have the following identity for number conservation
\beq
\sum_{i=1}^N m_i = K.
\eeq
We shall show that when the total number of clusters is $K$, the conditional probability distribution function is given by
\beqq
p'(m_1,\dots,m_N | K)   =   \frac{1}{C_{N,K}} \frac{a(1) ^{m_1}\dots a(N)^{m_N}}{m_1!...m_N!},
\eeqq
where the normalization constant $C_{N,K}$ can be computed explicitly (see formula \ref{CNKaconstant}). We will use this formula to compute the statistical moments for the cluster distributions. Using a combinatorial approach, we present several explicit formula for the steady-state distributions of clusters. Some of the results presented here were previously announced without proofs in the short letter \cite{Hoze2012}.

The paper is organized as follows. In section \ref{section:CFPmodels}, we present the stochastic model of coagulation-fragmentation for $N$ independent particles, that we analyze by a Markov chain. We obtain explicit formula for the cluster configuration combinatorics using the partition of the integer $N$, for the distribution of particles in clusters. We derive the time evolution equations for the number of clusters. These equations represent a novel Markov chain, which allows us to determine the number of clusters at steady-state. By combining the expression for number of clusters with the distribution of cluster in configuration conditioned on the number of clusters, we obtain the distribution of the particles in clusters. In section \ref{section:times}, we introduce two characteristic times for studying the time distribution two particles are in a cluster. These two times characterize the dynamics of exchange of particles between clusters: the first one is the mean time that two particles spend together before they separate, and the second is the mean time that they spend separated before they meet again for the first time. The colocalization probability of two particles is defined as the fraction of time, the particles spend together. Sections \ref{section:example1}-\ref{section:example2}-\ref{section:example3} provide direct application of these results to specific CFP: we focus specifically on the case of constant coagulation and fragmentation kernels in section \ref{section:example1}, for which we obtain analytical expressions for the clusters distributions. We also obtain several formula when the cluster size is limited.
%

{In this article, to determine statistics of the cluster distribution, we alternatively study two different systems. First, we study the distribution of clusters using the integer partition of the total number of particles. We obtain the probability distribution of cluster configurations. In this process, we use the term \textit{coagulation} when two clusters of a given size coalesce and form a new cluster. We use the term \textit{fragmentation} to describe the separation of a cluster into two smaller ones. The coagulation and fragmentation kernels that we have chosen here allow us to perform another analysis: when the number of clusters is fixed, the overall rates of coagulation and fragmentation are independent of the configurations of the clusters. We will study aggregation-fragmentation, where the number of clusters is known.  In that case, we shall use the following terminology:  \textit{formation} describes the change when  a distribution of $K$ becomes $K-1$ clusters. We use \textit{separation} to describe the process by which a distribution of $K$ clusters is transformed into a distribution of $K+1$ clusters.}

%

\section{Coagulation-fragmentation with a finite number of independent particles} \label{section:CFPmodels}

\subsection{Stochastic coagulation-fragmentation equations for a finite number of particles}
To describe the steady-state distribution for a CFP stochastic model with a finite number of $N$ particles, we shall use a continuous-time Markov chain in the space of cluster configurations. The $N$ particles distributed in clusters of size {$(n_1,\dots,n_N)$} can undergo coagulation or fragmentation events under the constraint that
\beq
\sum_{k=1}^N n_k  =N.
\eeq
{ To study the distribution of particles in clusters, we use the decomposition of the integer $N$ as the sum of positive integers (integer partition) \cite{Andrews1976}. The partitions of the integer $N$ are described in N dimensions by the ensemble
\beq
P_N= \{(n_1,\dots,n_N) \in \mathbb{N}^N; \sum_{i=1}^N n_i=N   \mbox{ and } n_1 \geq \dots \geq n_N\geq 0\}.
\eeq}
The probability $P(n_1,\dots,n_N,t)$ of the configuration $(n_1,\dots, n_N)$ at time $t$ satisfies an ensemble of closed equations.  Indeed, by considering all the possible coagulation or fragmentation events, the Master equation is obtained by considering the events occurring between time $t$ and $t+ \Delta t$:
\begin{itemize}
\item Two clusters of size $n_i$ and $n_j$ coagulate with a probability  $C(n_i,n_j)\Delta t$ to form a cluster of size $n_i+n_j$.
\item A cluster of size $n_i$ dissociates into two clusters of size $k$ and $n_i-k$  with a probability $F(k,n_i-k)\Delta t$.
\item Nothing happens with the  probability {$1- \sum_{i=1}^{N-1} \sum_{j=i+1}^{N}C(n_i,n_j) \Delta t - \sum_{i=1}^{N} \sum_{k=1}^{n_i-1}F(k,n_i-k)\Delta t$}.
\end{itemize}
Thus, the Master equations are

{ \beq
\frac{d}{dt}P(n_1,\dots,n_N,t) &=& - \left(\sum_{i=1}^{N-1} \sum_{j=i+1}^{N}C(n_i,n_j) + \sum_{i=1}^{N} \sum_{k=1}^{n_i-1}F(k,n_i-k)\right)P(n_1,\dots,n_N,t) \nonumber \\
&+&\sum_{k=1}^N \sum_{\underset{n_i'+n_j'=n_k}{n_i'>0,n_j'>0}   } C(n_i',n_j')  P(n_1,\dots,n_i',\dots, n_j',\dots ,n_N,t) \nonumber \\
&+&\sum_{i=1}^{N-1} \sum_{j=i+1}^N F(n_i,n_j)   P(n_1,\dots,n_i+n_j,\dots ,n_N,t).
\eeq
}
{ Moreover,  $C(n_i,n_j)=0$ if either $n_i$ or $n_j$ is equal to 0. }
We now introduce an ensemble in $\mathbb{N}^N$ which consists of integer decompositions of clusters with a given size:
  \beq
{P_N'= \{(m_1,\dots,m_N) \in \mathbb{N}^N; \sum_{i=1}^N i m _i=N   \mbox{ and } m_1 , \dots , m_N\geq 0\}.
}\eeq
In the ensemble $P_N'$,  $m_i$ is the number of occurrence of the integer $i$ in the partition of the integer $N$.
{ The two ensembles $P_N$ and $P_N'$ correspond to different representations of the clusters distributions. For example $N=9$ particles are distributed in two clusters of one particle, two clusters of two, and one cluster of three and the distribution is $(3,2,2,1,1,0,0,0,0) \in P_9$, and $(2,2,1,0,0,0,0,0,0) \in P_9'$.}

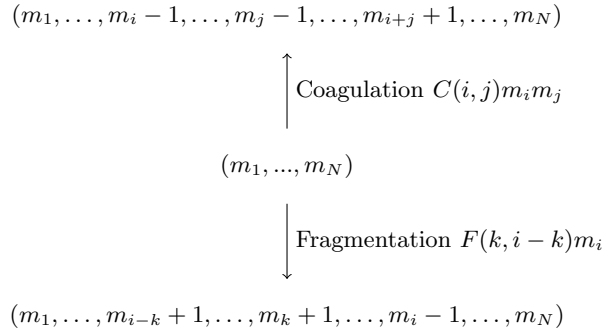
\begin{figure}[http!]
    \centering
\begin{tikzpicture}
\node at (0,0) {$(m_1,...,m_N)$};
\draw[->] (0,0.5) -- (0,1.5);
\node[right] at (0,1) {Coagulation $C(i,j)m_im_j$ };
\node at (0,2) {$(m_1,\dots,m_{i}-1,\dots,m_{j}-1,\dots,m_{i+j}+1,\dots,m_N)$};
\draw[->] (0,-0.5) -- (0,-1.5);
\node at (0,-2) {$(m_1,\dots,m_{i-k}+1,\dots,m_{k}+1,\dots,m_{i}-1,\dots,m_N)$};
\node[right] at (0,-1) {Fragmentation $F(k,i-k)m_i$ };
\end{tikzpicture}
\caption{{\bf Markov chain representation of the transitions}. Starting with a configuration $\mn$, $m_i$ is the number of clusters of size $i$ and $N$ is the total number of particles. The fragmentation rate of a cluster of size $i$ into one cluster of size $k$ and one of size $i-k$ is $F(k,i-k)m_i$, while the rate of formation of a cluster of size $i+j$ from two clusters of size $i$ and $j$ is equal to $C(i,j)m_im_j$. }\label{fig:Mathm1mN}
\end{figure}

 A  sufficient condition to obtain an invariant measure of the steady state probability is the reversibility of the CFP \cite{Kelly1979,Liggett1985} where the coagulation-fragmentation kernel satisfies the detailed balance condition: there exists a function $a(i)=a_i$ such that \cite{Durrett1999}
\beq\label{aii}
C(i,j)a_ia_j = F(i,j)a_{i+j}.
\eeq
{ The functions $a_i$ characterize the ratio of the coagulation and fragmentation rates,  and any function of the form $a'_i = \alpha^i a_i$ with $\alpha \neq 0$ satisfies the above condition. }
This condition insures the reversibility of the Markov chain and guarantees the existence of an invariant measure \cite{Durrett1999}, where the steady-state probability of a given configuration $\mn \in P_{N}'$ is
\beq
P' \mn= \frac{1}{C_N} \frac{a_1^{m_1}...a_N^{m_N}}{m_1!...m_N!},
\label{Probadistribution}
\eeq
where $C_N$ is a normalization constant. An explicit computation of the  { normalization constant}  is difficult \cite{Thompson1988}. Here, we propose to estimate the probability of occurrence of a certain cluster configuration $\mn$ by limiting the study to the configurations of a given number of particles.\\
At this stage, we shall explain the rational for computing expression \eqref{Probadistribution}. This expression is the distribution at equilibrium of particles in clusters, where the dissociation (resp. association) rate is proportional to the number of elements (minus one) (resp. the number of pairs of particles). \\
The equilibrium probability distributions associated to the Markov chain configuration $\mn$ is computed from analyzing the transition between the two neighboring states $(m_1,\dots,m_i-1,\dots,m_j-1, \dots,m_{i+j}+1,\dots,m_N)$ and $\mn$. { It is obtained first from the coagulation rate  $\psi(i,j) $ of a cluster of size $i$ with one of size $j$, given the distribution $\mn$. The rate $\psi(i,j) $ is given by}
\beq
\psi(i,j) &=& \frac{1}{2}C(i,j)m_im_j\, \hbox{ if } i\neq j \\
          &=& C(i,i)m_i(m_i-1)\,  \hbox{ otherwise. }
\eeq
The factor $\frac{1}{2}$ accounts for the   {symmetric} cases $\psi(i,j)$ and $\psi(j,i)$. {The fragmentation rate $\phi(i,j)$ from $i+j$ to $(i,j)$,  that accounts for the transition from the configuration $(m_1,\dots,m_i-1,\dots,m_j-1, \dots,m_{i+j}+1,\dots,m_N)$ to $\mn$ (Fig. \ref{fig:Mathm1mN}) is}
\beq
\phi(i,j) = F(i,j) (m_{i+j}+1).
\eeq
Thus, the stationary probability $\pi'$ satisfies the relation
\beq
 \frac{\pi'(m_1, \dots,m_i-1,m_j-1,m_{i+j}+1, \dots,m_N)}{\pi' \mn} &=& \frac{\psi(i,j)}{\phi(i,j)} \nonumber \\
& =&\frac{1}{2}\frac{C(i,j)}{F(i,j)}\frac{m_im_j}{m_{i+j}+1}  \nonumber \\
& =&\frac{1}{2}\frac{a_ia_j}{a_{i+j}}\frac{m_im_j}{m_{i+j}+1} .
 \label{eq:Mathprobabilitypi}
 \eeq
A direct computation shows that the probability $P' \mn$ defined in eq.  \eqref{Probadistribution} satisfies equation \eqref{eq:Mathprobabilitypi}.

\newpage

\subsection{Cluster partitions with a finite number of particles}
 To determine the cluster distribution at equilibrium, we compute here the probability of a configuration when the number of clusters $K$ is fixed. We also  find the probability of having $K$ clusters. The number of distributions of $N$ particles into $K$ clusters is the
cardinal of the ensemble
\beq
\ANK= \{(n_1,\dots,n_K) \in (\mathbb{N}^*)^K; \sum_{i=1}^K n_i=N   \mbox{ and } n_1 \geq \dots \geq n_K\geq 1\},\nonumber \\
\eeq
which is also the ensemble of the partitions of the integer $N$ as a sum of $K$ integers. This ensemble is in bijection with
\beq
\ANK'= \{\mn\in \mathbb{N}^N; \sum_{i=1}^N im_i=N \mbox{ and } \sum_{i=1}^{N} m_i=K\},
\eeq
where the application $\ANK \rightarrow \ANK'$   defined by
{\beq
(n_1,\dots,n_K)  \mapsto  (m_1,\dots, m_N) = \large( \sum_{i=1}^K 1_{\{n_i=1\}},\dots,  \sum_{i=1}^K 1_{\{n_i=N\}}\large)
\eeq}
{maps the partition $(n_1,\dots,n_K)$ where $N$ is written as a sum of $K$ positive  integers to the number of occurrence of each integer into the image partition. } The partitions of $N$ are written as
\beq
P_N = \bigcup_{K} \ANK \mbox{  and } P_N' =  \bigcup_{K} \ANK'.
\eeq
In section \ref{section:example1}, \ref{section:example2} and \ref{section:example3}, we derive explicitly expressions for the probabilities of configurations in $\ANK'$.
\subsection{Statistical moments for the cluster configurations when the number of clusters is fixed}
We show now that the  probability of configuration $\mn$, when the total number of clusters is equal to $K$, is given by
\beq
 p'(m_1,\dots,m_N|K) =\frac{\ds{\frac{a_1^{m_1}...a_N^{m_N}}{m_1!...m_N!}}}{C_{N,K}},
\label{eq:Mathfproba}
\eeq
where
\beq\label{cnk}
C_{N,K}=\sum_{(m_i) \in \ANK'} \ds{\frac{a_1^{m_1}...a_N^{m_N}}{m_1!...m_N!}}.
\eeq
The normalization factor of eq. \eqref{eq:Mathfproba} is computed using the following result:
\begin{remark} \label{remarkgeneratingfunction}
We consider the functions
\beq
S(x) = \sum_{i=1}^{\infty}a_ix^i
\label{generatingfunction}
\eeq
and the partial sums
\beq
S_N(x) = \sum_{i=1}^{N}a_ix^i.
\eeq
The $K^{th}-$power of these functions, that is $S^K$ and $S_N^K$ have the same $N^{th}$ order coefficient and this coefficient determines $C_{N,K}$. Moreover, the function
\beq
g(x,y)=\exp(S(x)y)
\label{generatingexpS}
\eeq
is a generating function of the $C_{N,K}$.
\end{remark}
\Proof
The number of configurations $\mn$ that satisfy the conditions $\sum_i m_i=K$ and $\sum_i im_i=N$ is the $N^{th}$-order coefficient of the multinomial expansion
\beqq
(a_1x_1+\dots+a_Nx_N)^K=\sum_{\mn, \,\sum m_i=K} \frac{K!}{m_1!\dots
m_N!}(a_1x_1)^{m_1}\dots (a_Nx_N)^{m_N}.
\eeqq
Using the $N$-tuple $(x_1,x_2,\dots, x_N)=(X,X^2,\dots,X^N)$, we obtain
\beq \label{sum1}
\frac{(a_1X+\dots+a_NX^N)^K}{K!}=\sum_{\mn, \,\sum m_i=K}\prod_i{  \frac{ a_i^{m_i}}{m_i!} }X^{\sum_i im_i}.
\eeq
We can group the terms by the exponents of $X$, which are equal to $\sum_i im_i$.  In particular, for a partition $\mn \in \ANK'$, the exponent is $\sum_i im_i = N$,  and the $N-$th order coefficient in expression \eqref{sum1} is equal to
\beq  \sum_{\mn \in \ANK'} \frac{\prod a_i^{m_i}}{m_1!\dots
m_N!} =C_{N,K}.
\eeq
%
More generally,
\beq
\frac{S^K(X)}{K!} &=&  \sum_{n=K}^{\infty} \sum_{\mn \in P_{n,K}'} \frac{ \prod a_i^{m_i}}{m_1!\dots
m_N!}X^{n} \\
&=&  \sum_{n=K}^{\infty} \sum_{\mn \in P_{n,K}'}C_{n,K}X^{n} .
\label{equationSK}
\eeq
It follows that the function defined by
\beq
g(X,Y) &=& \exp(S(X)Y) \\
&=& \sum_{K=0}^{\infty} \frac{S^K(X)}{K!}  Y^K\\
&=& \sum_{K=0}^{\infty}  \sum_{n=K}^{\infty}C_{n,K}X^{n}Y^K
\eeq
is  a generating function of the $C_{N,K}$. \\

\begin{remark}
The coefficients $C_{N,K}$ defined by relation \eqref{cnk} satisfy the induction formula:
\beq
(N+1)C_{N+1,K} = \sum_{k=0}^{N-K+1}(k+1) a_{k+1}C_{N-k,K-1}.
\label{recurrentrelationCNK}
\eeq
with
\beq
\left \lbrace
\begin{array}{ccc}
C_{N,N} &=&\ds{ \frac{a_1^N}{N!}} \\ &&\\
C_{N,1} &=& a_N.
\end{array}\right.
\label{CNNCN1}
\eeq
\end{remark}

\Proof
We start with the formulas \eqref{CNNCN1}. The coefficient $C_{N,N} $ is obtained from the unique partition in $P_{N,N}'$, which is given by $m_1=N$ and $m_i=0$ if $i>1$. Therefore, $C_{N,N}  =\frac{a_1^N}{N!} $. The coefficient $C_{N,1} $ is obtained from the partition  $m_N=§$ and $m_i=0$ if $i<N$, and thus $C_{N,1} =  a_N$.

We now prove relation \eqref{recurrentrelationCNK}.  Differentiating the function $\sigma(x)  = \frac{S^K(x)}{K!}  $, we get
\beq
\sigma'(x)   = S'(x)  \frac{S^{K-1}(x)}{(K-1)!}.
\label{sigmaprime}
\eeq
We evaluate the l.h.s of \eqref{sigmaprime} using  eq. \eqref{equationSK}, and obtain
 \beq
 \sigma'(x)  &=& \left( \sum_{n=K}^{\infty} C_{n,K}x^n \right)' \nonumber \\
 &=& \sum_{n=K-1}^{\infty}  (n+1) C_{n+1,K}x^{n} \nonumber \\
  &=&x^{K-1} \sum_{n=0}^{\infty}  (n+K) C_{n+K,K}x^{n}.
 \eeq
We evaluate the r.h.s of \eqref{sigmaprime} using the definition of $S$ \eqref{generatingfunction}  and we obtain
\beq
 \sigma'(x)   &=& \left( \sum_{i=0}^{\infty}(i+1) a_{i+1} x^i  \right)  \left( \sum_{i=K-1}^{\infty} C_{i,K-1}x^i \right) \nonumber \\
 &=& x^{K-1} \sum_{n=0}^{\infty}  \left( \sum_{k=0}^{n}(k+1) a_{k+1} C_{n-k+K-1,K-1}  \right)x^n .
\eeq
Thus, by equalizing the $N-$th order coefficient of $\sigma'(x)$, we have
\beq
(N+1)C_{N+1,K} = \sum_{k=0}^{N-K+1}(k+1) a_{k+1}C_{N-k,K-1}.  \mbox{ \Endproof }
\eeq
Next, we estimate various moments when the number of clusters is fixed. We summarize the main result in the
\begin{thm} \label{thmClusterSize}
When the number of clusters is equal to $K$ for a total of $N$ particles, the mean number of clusters of size $i$ is
\beq
\langle M_i \rangle_{N,K} =  a_i \frac{C_{N-i,K-1}}{C_{N,K}},
\label{eq:MiNK}
\eeq
where $a_i$ and $C_{N,K}$ are defined in \eqref{aii} and \eqref{cnk} respectively.
\end{thm}
\Proof
The mean number of clusters of size $i$ when the total number of clusters is $K$ is given by the following sum
\beq
\langle M_i \rangle_{N,K}&=& \sum_{\ANK'}m_ip'(m_1,...,m_N) \\
&= & \frac{1}{C_{N,K}}\sum_{\ANK'}m_i\frac{a_1^{m_1}...a_N^{m_N}}{m_1!...m_N!}\nonumber  \\
&=&  \frac{1}{C_{N,K}}\sum_{\ANK' , m_i>0}a_i\frac{a_1^{m_1}...a_i^{m_i-1}...a_N^{m_N}}{m_1!...(m_i-1)!...m_N!},
\eeq
where the subscript $\ANK' , m_i>0$ in the sum characterizes the partitions in $\ANK'$ containing  at least one occurrence of the integer $i$.
By considering the partitions of $\ANK'$ where $i$ appears at least once and removing one $i$, we obtain exactly the partitions of $P_{N-i,K-1}$, except for the number of occurrence of $i$, where the corresponding partitions in both
sets have the same number of repetitions, i.e.
\beq
\forall j\neq i, m_j \hbox{ in } (m_1,\dots,m_{N-i}) \in P'_{N-i,K-1}
\eeq
is equal to
\beq
m_{j} \hbox{ in } (m_1,\dots,m_{N}) \in  \ANK'
 \eeq
and $m_i$ in $(m_1,\dots,m_{N-i}) \in P'_{N-i,K-1}$ is equal to $m_i+1$ in $(m_1,\dots,m_N) \in \ANK'$. There are no clusters larger than $N-i$ in the partitions ($m_j=0$ for $j>N-i$ for  $(m_1,\dots,m_N) \in \ANK' , m_i>0$). Thus
\beq
\langle M_i \rangle_{N,K}&=&  \frac{1}{C_{N,K}}\sum_{\ANK' , m_i>0}a_i\frac{a_1^{m_1}...a_i^{m_i-1}...a_N^{m_N}}{m_1!...(m_i-1)!...m_N!} \nonumber \\
&=&  \frac{a_i}{C_{N,K}}\sum_{ P_{N-i,K-1} } \frac{a_1^{m_1}...a_i^{m_i}...a_N^{m_N}}{m_1!... m_i!...m_N!}  \nonumber \\
&=&   a_i \frac{C_{N-i,K-1}}{C_{N,K}} . \mbox{ \Endproof }
\label{Migeneralcase}
\eeq
We further have:
\begin{remark} \label{remarkmi0}
$\langle M_i \rangle_{N,K}=0$ if $i>N-K+1$.
\end{remark}
\Proof
When the $N$ particles are distributed in $K$ clusters, the largest cluster contains at most $N-K+1$ particles. The corresponding partition in $\ANK'$ is given by $m_1=K-1$, $m_{N-K+1}=1$, and $m_j=0$ otherwise.
\Endproof

\begin{thm}
The second moment of the number of clusters of size $i$ is
\beq
\langle M ^2_i\rangle_{N,K} &=& \frac{1}{C_{N,K}}\sum_{\ANK'}m^2_i\frac{a_1^{m_1}...a_N^{m_N}}{m_1!...m_N!} \\
&=&   a^2_i \frac{C_{N-2i,K-2}}{C_{N,K}}  + a_i \frac{C_{N-i,K-1}}{C_{N,K}}, \nonumber
\eeq
and the covariance is
\beq
\langle M ^2_{i,j}\rangle_{N,K} - \langle M ^2_{j}\rangle_{N,K}\langle M _{j}\rangle_{N,K}=a_ia_j\left(\frac{C_{N-i-j,K-2}}{C_{N,K}} -\frac{C_{N-i,K-1}C_{N-j,K-1}}{C_{N,K}^2} \right) . \nonumber \\
\eeq
\end{thm}
\Proof
The variance of the number of clusters of size $i$ is given by
\beq
\langle M ^2_i\rangle_{N,K} &=& \frac{1}{C_{N,K}}\sum_{m_i \in \ANK'}m^2_i\frac{a_1^{m_1}...a_N^{m_N}}{m_1!...m_N!} \nonumber  \\
&=&     \frac{1}{C_{N,K}}\sum_{ m_i \in \ANK'} [m_i(m_i-1)+m_i ]\frac{a_1^{m_1}...a_N^{m_N}}{m_1!...m_N!} \nonumber \\
&=&     \frac{1}{C_{N,K}}\sum_{ m_i \in \ANK', m_i>1}a_i^2 \frac{a_1^{m_1} \dots a_i^{m_i-2} \dots a_N^{m_N}}{m_1!\dots (m_i-2)! \dots m_N!}  +\langle M_i \rangle_{N,K}.  \nonumber
 \eeq
Using an argument similar to the proof of Theorem \ref{thmClusterSize}, we obtain
\beq
\langle M ^2_i\rangle_{N,K} &=&     \frac{1}{C_{N,K}}\sum_{ m_i \in  P_{N-2i,K-2}}a_i^2 \frac{a_1^{m_1} \dots a_i^{m_i} \dots a_N^{m_N}}{m_1!\dots m_i   ! \dots m_N!}  +\langle M_i \rangle_{N,K} \nonumber \\ &&  \nonumber \\
&=&     a_i^2 \frac{C_{N-2i,K-2}}{C_{N,K}} +a_i \frac{C_{N-i,K-1}}{C_{N,K}}. \nonumber \\
\eeq
The covariance of the number of clusters of   size $i$ and $j$  is obtained from the term
\beq
\langle M ^2_{i,j}\rangle_{N,K} =\frac{1}{C_{N,K}}\sum_{\ANK'}m_i m_j\frac{a_1^{m_1}...a_N^{m_N}}{m_1!...m_N!},
\eeq
which, by the same reasoning, leads to
\beq
\langle M ^2_{i,j}\rangle_{N,K} =a_ia_j\frac{C_{N-i-j,K-2}}{C_{N,K}}.
\eeq
Using eq. \eqref{Migeneralcase},  the covariance is
\beq
\langle M ^2_{i,j}\rangle_{N,K} - \langle M ^2_{j}\rangle_{N,K}\langle M _{j}\rangle_{N,K}=a_ia_j\left(\frac{C_{N-i-j,K-2}}{C_{N,K}} -\frac{C_{N-i,K-1}C_{N-j,K-1}}{C_{N,K}^2} \right) .  \nonumber \\
\eeq
\Endproof
%
%
\section{Distribution of the number of clusters}
\begin{figure}[http!]
\begin{center}
\includegraphics[scale=0.6]{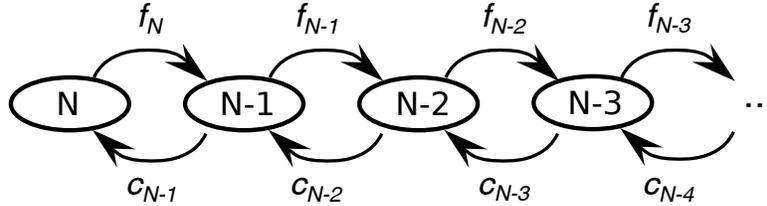}
\caption{ {\bf Markov chain representation for the number of clusters.} $s_K$ (respectively $f_K$) is the separation (respectively formation) rate of a cluster when there are $K$ clusters.}
\label{fig:schemageneral}
\end{center}
\end{figure}

In the previous section, we computed the probability distribution of a cluster configuration and determined the statistical moments for a fix number of clusters. In this section, we study the statistics of the entire cluster configurations. We focus here on the probability distribution of the \textit{number of clusters} and we shall compute the time dependent probability density function
\beq
P_{K}(t)= P\{ \mbox{$K$ clusters at time } t \},
\eeq
which is a birth-and-death process that we investigate using a Markov chain. \\
The probability of having $K$ clusters at time $t+\Delta t$ is the sum of the probability of starting at time $t$ with $K-1$ clusters and one of them dissociates into two smaller ones plus the probability of starting with $K+1$ clusters and two of them associate plus the probability of starting with $K$ and nothing happens (Fig. \ref{fig:schemageneral}).
The first probability is the product of $P_{K-1}$ by the transition rate $s_{K-1}\Delta t$ to go from state with $K-1$ clusters to $K$, while the second is the transition from $K+1$ to $K$, which is the product of $P_{K+1}$ by the transition rate $f_{K+1}\Delta t$ of going from $K+1$ clusters to $K$. Thus the master equations are given by
\beq
\label{eq:MathPKt}
\left \lbrace
\begin{array}{ccc}
\dot{P}_1(t) &=&-s_1P_1(t)+f_2P_{2}(t)\\ &&\\
\dot{P}_K(t)
&=&-(f_{K}+s_K)P_K(t)+f_{K+1}P_{K+1}(t)+s_{K-1}P_{K-1}(t)\\
&& \label{eq:MathMarkov1}\\
\dot{P}_N(t) &=&-f_NP_N(t)+s_{N-1}P_{N-1}(t).
\end{array}\right.
\eeq
%
In section \ref{section:example1},  \ref{section:example2} and \ref{section:example3},  we will solve this system of equations explicitly at steady-state for particular formation and separation kernels. We now derive a general formula for the steady probability
\beq
\Pi_K= \underset{ t\rightarrow \infty}{\lim}P_K(t)
\eeq
of having $K$ clusters at steady state and express it in terms of the $a_i$. The steady state probabilities of the number of clusters  are solution of the system
\beq
 \left \lbrace
\begin{array}{ccc}
0&=&-s_1 \Pi_1+f_2\Pi_{2}\\ &&\\
0&=&-(f_{K}+ s_K)\Pi_K+f_{K+1}\Pi_{K+1}+s_{K-1}\Pi_{K-1}\\
&& \label{eq:MathPiKEquilibrium}\\
0&=&-f_N \Pi_N +s_{N-1} \Pi_{N-1},
\end{array}\right.
\eeq
with the normalization condition
\beq
\sum_{K=1}^N\Pi_K=1.
\eeq
The probabilities $\Pi_K$ are given by the ratio
\beq
\frac{\Pi_{K} }{\Pi_{K-1}} = \frac{s_{K-1}}{f_{K}}  \mbox{ for } K\geq 2
\eeq
and the coefficients $s_K$ and $f_K$ are the mean-field separation and formation rates respectively. Whereas the cluster configurations when the number of clusters is fixed depend only on the kernel $a_i$, the statistics of the number of clusters depend on the cluster fragmentation and coagulation rates $F$ and $C$. In the following, we will focus on the coagulation condition $C(i,j)=1$ and the fragmentation $F(i,j) =\frac{a_ia_j}{a_{i+j}}$ to state the
\begin{thm}
When $C(i,j)=1$ and $F(i,j) =\frac{a_ia_j}{a_{i+j}}$, the separation rate when there are $K$ clusters is given by
\beq
s_K =  \frac{\sum_{i=1}^N \sum_{j=1}^{i-1}a_ja_{i-j} C_{N-i,K-1}}{C_{N,K}}
\eeq
and the formation rate when there are $K$ clusters is
\beq
f_K = \frac{K(K-1)}{2}.
\eeq
\end{thm}
\Proof
The total dissociation rate $d(n)$ of a cluster of size $n$ is obtained by summing over all the possible sizes resulting from the dissociation   and is given by
\beq
d(n) &=& \sum_{i=1}^{n-1}F(i,n-i) \nonumber \\
&=& \sum_{i=1}^{n-1} \frac{a_i a_{n-i}}{a_n}.
\label{dissociationrate}
\eeq
The  rate at which a given configuration $\mn \in \ANK'$ dissociates is $\sum_{i=1}^Nd(i)m_i$. The separation rate for $K$ clusters is thus
\beq
s_K&=&\sum_{\ANK'} \sum_{i=1}^Nd(i)m_ip'(m_1,...,m_N) \nonumber \\
&=&\sum_{i=1}^Nd(i) \langle M_i \rangle_{N,K} \nonumber \\
&=& \frac{\sum_{i=1}^Nd(i)a_i C_{N-i,K-1}}{C_{N,K}}.
\label{dissociationrateK}
\eeq
The separation rate of a distribution of $K$ clusters eq. \eqref{dissociationrateK} can thus be expressed as a function of the $C_{N,K}$
\beq
s_K =  \frac{\sum_{i=1}^N \sum_{j=1}^{i-1}a_ja_{i-j} C_{N-i,K-1}}{C_{N,K}}.
\eeq
For $K=1$, the only cluster is of size $N$ and the separation rate is
\beq
s_1  =d(N).
\eeq
The formation rates are given by
\beq
f_K= \sum_{\ANK'} \frac{1}{2}\left(\sum_{i=1}^{N}m_i(m_i-1) C(i,i)+\sum_{i\neq j}m_i m_j C(i,j) \right)   p'(m_1,...,m_N|K).  \nonumber \\
\eeq
For a distribution of $K$ clusters, this is equal to the number of cluster pairs
\beq
f_K =\frac{K(K-1)}{2}. \mbox{ \Endproof}
\eeq
We are now in position to study the statistics of the entire cluster configurations. Indeed, using Bayes' rule, the probability of a configuration $\mn$, that contains $K$ clusters is the product of the conditional probability $p'(m_1,\dots,m_N |K)$ by the probability of having $K$ clusters
\beq
p'(m_1,\dots,m_N,K )=p'(m_1,\dots,m_N |K)\Pi_K.
\eeq
The mean number of clusters of size $i$ is thus
\beq
\langle M_i \rangle_N = \sum_{K=1}^N \Pi_K \langle M_i \rangle_{N,K}.
\eeq
\section{Invariant of   clusters dynamics} \label{section:times}
\begin{figure}[http!]
\begin{center}
\includegraphics[scale=0.7]{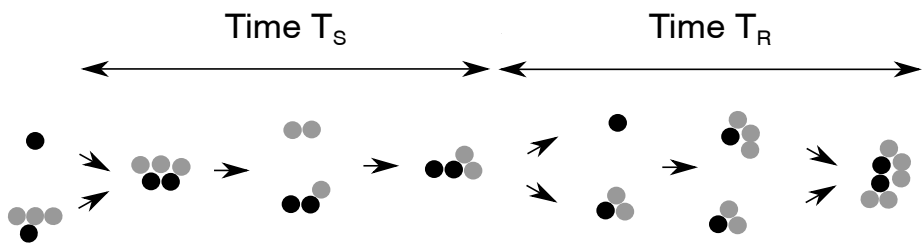}
\caption{ \textbf{Cluster formation and dissociation.} The time to separation $T_S$ is the mean time for two specific particles (black) to spend in the same cluster. Before physical separation, the cluster can aggregate with other clusters (grey) or some particle except either of the two can be separated. The time to association $T_R$ is the mean time the two particles meet again after separation for the first time.}
\label{fig:timetogether}
\end{center}
\end{figure}

We introduce and compute here several measures of the cluster configurations { that appear when the system is in a global steady state}. First, we compute the probability to find two particles in the same cluster, and second we measure two time scales associated to particle dynamics in an ensemble of clusters: 1) the mean time that two particles spend together before they separate and 2) the mean time that they spend separated before they meet again for the first time (Fig. \ref{fig:timetogether}). { The probability that two given particles are together is of interest in several cell biology examples: for instance, some genes can be silenced if the telomeres that carry them are forming a cluster \cite{Ruault2011}.}

\subsection{ The probability to find two particles in the same cluster}
When the mean number of clusters has reached its equilibrium, particles can still be exchanged between clusters. To characterize this exchange, we compute the probability to find two particles in the same cluster.  \\
When the distribution of the clusters is{ $(n_1,\dots,n_N)$}, the probability $P_2(n_1,\dots,n_N)$ to find two given particles in the same cluster is obtained by using the probability to choose the first particle in the cluster $n_i$, which is equal to the number of particles in the cluster divided by the total number of particles $\frac{n_i}{N}$. The probability to have the second particle in the same cluster is $\frac{n_i-1}{N-1}$. Summing over all possibilities, we get
\beq \label{Matheqq}
P_2(n_1,\dots,n_N)=\sum_{i=1}^K
\frac{n_i}{N}
\frac{n_i-1}{N-1}=\frac{1}{N(N-1)}(\sum_{i=1}^Kn^2_i-N).
\eeq
{ This probability is similar to Simpson's diversity index \cite{Simpson1949}, a measure frequently used to quantify the diversity of ecosystems.}
We note that
\beq
\sum_{j=1}^{N}n^2_j =  \sum_{i=1}^{N}i^2m_i.
\eeq
{Thus, when the distribution $(n_1,\dots, n_N)$ contains $K$ clusters, we use $(n_1, \dots,n_K)\in
\ANK$ and obtain by summing over all configurations of $K$ clusters}
\beq
\sum_{(n_1, \dots,n_K)\in
\ANK}p(n_1,\dots, n_K | K)\sum_{j=1}^{K}n^2_j&=&\sum_{(m_i)\in
\ANK'}p'(m_i)\sum_{j=1}^{N}j^2m_j \\
&=&  \sum_{j=1}^{N}j^2 \sum_{(m_i)\in
\ANK'}m_jp'(m_i)\\
&=& \sum_{j=1}^{N}j^2 \langle M_j \rangle_{N,K},
\label{P2Mj}
\eeq
where $\langle M_j \rangle_{N,K}$ is the mean number of clusters of size $j$, when there are $N$ particles distributed in $K$ clusters eq. \eqref{eq:MiNK}.
Taking into account all possible distributions of clusters, we obtain that the probability $\langle P_2 \rangle $ to find two particles in the same cluster is
\beq
\langle P_2 \rangle=\sum_{K=1}^N\sum_{(n_1,\dots,n_K)\in \ANK} P_2(n_1,\dots,n_K)p(n_i)\Pi_K,
\eeq
which can be written, using expressions \eqref{Matheqq} and \eqref{P2Mj} as
\beq
\langle P_2 \rangle=\frac{1}{N(N-1)}\sum_{K=1}^N\Pi_K \sum_{j=1}^{N}j^2 \langle M_j \rangle_{N,K}-\frac{1}{N-1}.
 \label{eq:P2generalcase}
\eeq
This approach can be generalized to the probability of having $n\geq 2$ particles together.

\subsection{Mean time for two particles to stay together in a  cluster}
\label{sec:meantime}
{ We define the mean time to separation (MTS) as the mean time between the arrival of two given particles in a cluster and their separation after dissociation of the cluster.}
We first compute the probability of a configuration $\mn$ conditioned on having two particles in the same cluster. We then derive the transition rates from this particular configuration to any of the states accessible by a single dissociation or association event (see Fig. \ref{fig:Mathm1mN}). The accessible states are divided into two classes: first, the configurations for which the particles are not initially in the same cluster (separated state) and second the ones where they are together.
{ The current state is the configuration $\mn$, for which the particles are in the same cluster. Upon a coagulation event the particles stay in the same cluster. }
{ A dissociation event can either occur for a cluster that did not contain the tracking particles or for the cluster that contained the particles. In the latter case, there are two possibilities:  either the particles stay together after dissociation, or they are redistributed to two different clusters (separated state).}
The rate of dissociation from the ensemble $\mn$ to $(m_1,\dots,m_i+1,\dots,m_{k-i}+1,\dots,m_{k}-1,\dots, m_N)$ is  $2F(i,k-i)m_k$ if $i\neq k/2$ and $F(\frac{k}{2},\frac{k}{2})m_k$ otherwise.

\begin{proposition}
The probability that two particles in a cluster of size $k$, { given a configuration
$\mn \in \ANK '$} separate after a dissociation is
\beq
p_S(k; m_1, \dots, m_N)=\frac{1}{k(k-1)}\frac{1}{\sum_{i=1}^Nd(i)m_i} \sum_{i=1}^{k-1}i(2k-i)F(i,k-i). \nonumber \\
\eeq
\end{proposition}
\Proof
The probability that two particles in a cluster of size $k$ separates { after any dissociation event} is equal to the product of the probability that the particles are in the cluster multiplied by the probability that this dissociation results in the effective separation of the particles. The first probability is obtained by considering { the  total dissociation rate  $d(n)$ of a cluster of  size $n$ (see eq. \eqref{dissociationrate}). The probability that the cluster containing the two particles effectively  dissociates is thus proportional to $d(k)$ and is normalized by the total dissociation rate for the cluster configuration $\mn$ : }
\beq
\frac{d(k) }{\sum_{i=1}^Nd(i) m_i}.
\eeq
The probability of an effective separation $P_{sep}$ is the complementary of the probability that the two particles stay together. When the two resulting clusters are of size $i$ and $k-i$, this probability is equal to
\beq
P_{sep}=1- \frac{i(i-1)+(k-i)(k-i-1)}{k(k-1)}.
\eeq
{ The probability that the two particles initially in a cluster of size $k$ will be separated is thus equal to
}\beq
p_S(k)&=&\frac{d(k) }{\sum_{i=1}^Kd(i)m_i} \sum_{i=1}^{k-1}\left(1- \frac{i(i-1)+(k-i)(k-i-1)}{k(k-1)}\right)\frac{F(i,k-i)}{d(k)} \nonumber\\
&=&\frac{1}{k(k-1)}\frac{1}{\sum_{i=1}^Nd(i)m_i} \sum_{i=1}^{k-1}F(i,k-i)i(2k-i).
\mbox{\Endproof}
\eeq
The transition probability from the state $\mn$ to the separated one  equals the sum of the probabilities $p_S(k)$ to be separated over all cluster sizes $k$,
\beq
P_S\mn=\sum_{k=1}^Kp_S(k)m_k.
\eeq
To derive the MTS, we first write the transition matrix of the Markov chain representing the transitions between the configurations of separated and non-separated states, and second, we   determine the mean transition times between the states.

Ordering the configurations $\mn$ by arbitrary indices $I$, we define the transition matrix $T$ of size $(q(N)+1) \times (q(N)+1)$, where $q(N)$ { is the total number of different partitions of the integer $N$}. The elements of $T$ of the first $q(N)$  rows and columns are the transition probabilities between states where the two particles are together, while index $q(N)+1$ represents the separated state.

For $I,J \in (1,\dots,q(N))$, the matrix entries $[T]_{I,J}$ are either the coagulation or dissociation rates given above (see Fig. \ref{fig:Mathm1mN}), while the transitions $[T]_{I,q(N)+1}$ are equal to the rates to the separated state. Because the separated state is absorbing, we finally set $[T]_{q(N)+1,q(N)+1}=1$. Additionally, the matrix is normalized into a stochastic matrix such that
\beq
\forall I, \sum_{J=1}^{q(N)+1} [T]_{I,J}=1.
\eeq
We now estimate the mean time the two particles stay together. The mean time $\tau\mn$ from a configuration $\mn$ to any configuration $(m_1',\dots,m_N')$ accessible by a single coagulation or fragmentation event is {the reciprocal of the sum of all transition rates}
\beq
\tau\mn  = \bigg( \sum_{k=1}^{N}d(k)m_k+ \frac{K(K-1)}{2} \bigg)^{-1},
\eeq
{where $K = \sum_{i=1}^{N}m_i$ and the coagulation kernel is constant $C(i,j)=1$.
We represent the transition times in a vector  $\boldsymbol{\tau} $
{ \beq
\boldsymbol{\tau}  =   \begin{bmatrix}
    \tau(N,0,\dots,0)              \\[0.3em]
    \tau(N-2,1,\dots,0)       \\[0.3em]
    \vdots     \\[0.3em]
       \tau(0,\dots,0,1)
     \end{bmatrix}.
 \eeq}
We shall now compute the vector
{ \beq
  \boldsymbol{t}  =   \begin{bmatrix}
    t(N,0,\dots,0)              \\[0.3em]
    t(N-2,1,\dots,0)       \\[0.3em]
    \vdots     \\[0.3em]
       t(0,\dots,0,1)
     \end{bmatrix}
 \eeq}
 which is the MTS for an ensemble of cluster configuration $\mn$: it is related to the vector $\boldsymbol{\tau}$ by \cite{Meyer1975}
\beq
\boldsymbol{t}  =\sum_{n=0}^{\infty}T^n\boldsymbol{\tau}.
\label{MFTS}
\eeq
%
The vector $ \boldsymbol{t}  $ is computed using the matrices
\beq
A=I_{q(N)+1} -T
\eeq
and
\beq
A^* =\left[A_{q(N)+1}\right]^{-1},
\eeq
where $ A_{q(N)+1}$ is the matrix $A$ from which the $q(N)+1$ row and column were removed. Equation \eqref{MFTS} is equivalent to
}
\beq
\boldsymbol{t}= A^*\boldsymbol{\tau}.
\eeq
The MTS averaged over the equilibrium configuration distribution is
\beq
T_S = \boldsymbol{p}^{*T}\boldsymbol{t}
\eeq
where $\boldsymbol{p}^{*T} = [ p^*\mn]_{\mn}$ and $p^*\mn$ is the probability distribution of the configuration $\mn$ when the two particles are in the same cluster. Using  Bayes theorem, the probabilities  $p^*\mn$ are given by
\beq
p^*\mn  =   \frac{  P_2\mn p\mn}{\langle P_2 \rangle},
\eeq
where we recall that $p\mn$ is the probability of the configuration $\mn$, $\langle P_2 \rangle$ is the probability that two specific particles are in the same cluster (see computation eq. \eqref{eq:P2generalcase}), and $P_2\mn$ is the probability that the two particles are in the same cluster configuration $\mn$. { The conditional probability  $p^* (m_1,\dots,m_N)$ to select two particles in the same cluster, conditioned on the distribution $(m_1,\dots,m_N)$ is larger for clusters of larger sizes}.
We conclude this section with the following
\begin{remark}
The time that two particles spend separated $T_R$ can be similarly determined, only the absorbing state in the transition matrix is the state at which the two particles are in the same cluster. Interestingly, the probability two particles are together is the fraction of time they spend together and is given by the ratio
\beq
\langle P_2 \rangle= \frac{T_S}{T_S+T_R}.
\eeq
\end{remark}
In the rest of the manuscript, we apply the previous analysis to three examples of coagulation-fragmentation with a finite number of particles.
\section{Example 1: the case $a_i=a$} \label{section:example1}
We consider the case of a constant kernel $a_i=a$. To compute the separation and formation rates $s_K$ and $f_K$, we use that $F(i,j)=a$ and $C(i,j)=1$. This fragmentation kernel corresponds to the following model: a cluster of size $n$ dissociates at a rate $\sum_{i=1}^{n-1}F(i,n-i)= (n-1)a$ and the sizes of the resulting clusters are uniformly distributed between 1 and $n-1$. {When there are $N$ particles and a total number of clusters $K$, the partition of clusters is denoted $(n_1,\dots,n_K) \in P_{N,K}$}. The total transition rate from a configuration of $K$ to $K+1$ clusters is the sum over all possible dissociation rates
 \beq
s_K =  \sum_{i=1}^{K}(n_i-1)a=(N-K)a.
\eeq
The formation rate is proportional to the number of pairs
\beq
{f_K= \sum_{i=1}^{K-1} \sum_{j=i+1}^{K}C(n_i,n_j) = \frac{K(K-1)}{2}.
}\eeq
Following relation \eqref{eq:MathPKt}, the steady-state probability $\Pi_K$ for the number of clusters of size $K$ satisfies the time independent master equation
\beq
\left \lbrace
\begin{array}{ccc}
s_1\Pi_1&=& f_2\Pi_{2},\\ &&\\\mu_1
(f_{K}+s_K)\Pi_K &=& f_{K+1}\Pi_{K+1}+s_{K-1}\Pi_{K-1},\\
&&\\
f_N\Pi_N &=& s_{N-1}\Pi_{N-1},
\end{array}\right.
\eeq
which leads to relation
\beq
\Pi_{K+1}=(2a)^K\frac{(N-1)!}{K!(K+1)!(N-K-1)!} \Pi_1.
\label{eq:MathPik}
\eeq
Using the normalization condition $\sum_K\Pi_K=1$, the probability $\Pi_1$ can be expressed as a hypergeometric series
\beq
\Pi_1=\frac{1}{ _1F_1(-N+1;2;-2a)},
\label{probability1cluster}
\eeq
where
\beq
_1F_1(a;b;z) = \sum_{n=0}^\infty \frac{(a)_n}{(b)_n} \frac {z^n} {n!},
\eeq
is Kummer's confluent hypergeometric function (\cite{Abramowitz} pp. 503--535) and
\beq
(x)_n=x(x+1)...(x+n-1)
\eeq
is the Pochhammer symbol. The average number of clusters at steady state
\beq
\mu_1(a)&=& \sum_{K=1}^N K\Pi_K \nonumber \\
 &=& \Pi_1\frac{d}{dz}\bigg(z _1F_1(-N+1;2;z)\bigg)_{|z=-2a}.
\label{deriv1}
\eeq
The derivative of the Kummer's function is
\beq
\frac{d}{dz} \ _1F_1(a;b;z)=\frac{a}{b} \ _1F_1(a+1;b+1;z).
\label{derive1F1}
\eeq
Finally the mean number of clusters is expressed as
\beq
\mu_1(a)&=&1+a(N-1)\frac{_1F_1(-N+2;3;-2a)}{ _1F_1(-N+1;2;-2a)},\nonumber \\
&=&1+a(N-1)G_1,
\label{NumberClusterModel2}
\eeq
where we note $G_1$ the function  defined by
\beq
G_1 = \frac{_1F_1(-N+2;3;-2a)}{_1F_1(-N+1;2;-2a)}.
\label{definitionG1}
\eeq
More generally, we introduce the functions $G_i$ defined by
\beq
G_i = \frac{_1F_1(-N+1+i;2+i;-2a)}{_1F_1(-N+1;2;-2a)}.
\label{definitionGi}
\eeq
Following the procedure presented above, all moments of the probability distribution $\Pi_K$ can be computed and the $n^{th}$-order moment $\mu_n$ is expressed using  the operator $H$ defined by \beq
H(f)(z) = \frac{d}{dz}zf(z),
\eeq
by
\beq
\mu_n = \sum_{n=1}^N K^n\Pi_K=\frac{{H^{(n)}({_1F_1}(-N+1;2;z))}_{|z=-2a}}{_1F_1(-N+1;2;-2a)}.
\eeq
Using the differentiation formula for the hypergeometric function \eqref{derive1F1}, the moments $\mu_n$ can be written as
\beq
 \mu_n&=& \sum_{k=0}^{n}\alpha_k^n\frac{(N-1)!}{(k+1)!(N-1-k)!} 2^ka^kG_k,\\
 &=&  \sum_{k=0}^{n}\alpha_k^n \frac{\Pi_{k+1}}{\Pi_1}G_k,
\eeq
 where the coefficients $\alpha_k^n$ are given by
\beqq
\alpha_k^n = \left \lbrace
\begin{array}{ccc}
&k!\sum_{j=0}^{k/2} (-1)^j\ds{\frac{(k+1-j)^n+(j+1)^n}{(k-j)!}}& \mbox{ if $k$ is even,}\\
\\
& k!\sum_{j=0}^{(k-1)/2} (-1)^j\ds{\frac{(k+1-j)^n-(j+1)^n}{(k-j)!}}& \mbox{ if $k$ is odd,}\\
\end{array}\right.
\eeqq
and $\alpha_0^n= \alpha_n^n=1$.
We can thus obtain the variance of the number of clusters, given by
\beq
\langle V_{\infty}(a) \rangle&=& \mu_2-\mu_1^2  \nonumber
\\
&=& a(N-1)G_1(a,N) \nonumber \\
&+&  \frac{2}{3}a^2(N-1)(N-2)G_2(a,N) \nonumber \\
&-& a^2(N-1)^2 G_1^2(a,N).
\eeq
\subsection{Asymptotic formulas for the mean and variance of the cluster number}
\begin{figure}[http!]
\begin{center}
\includegraphics[scale=1]{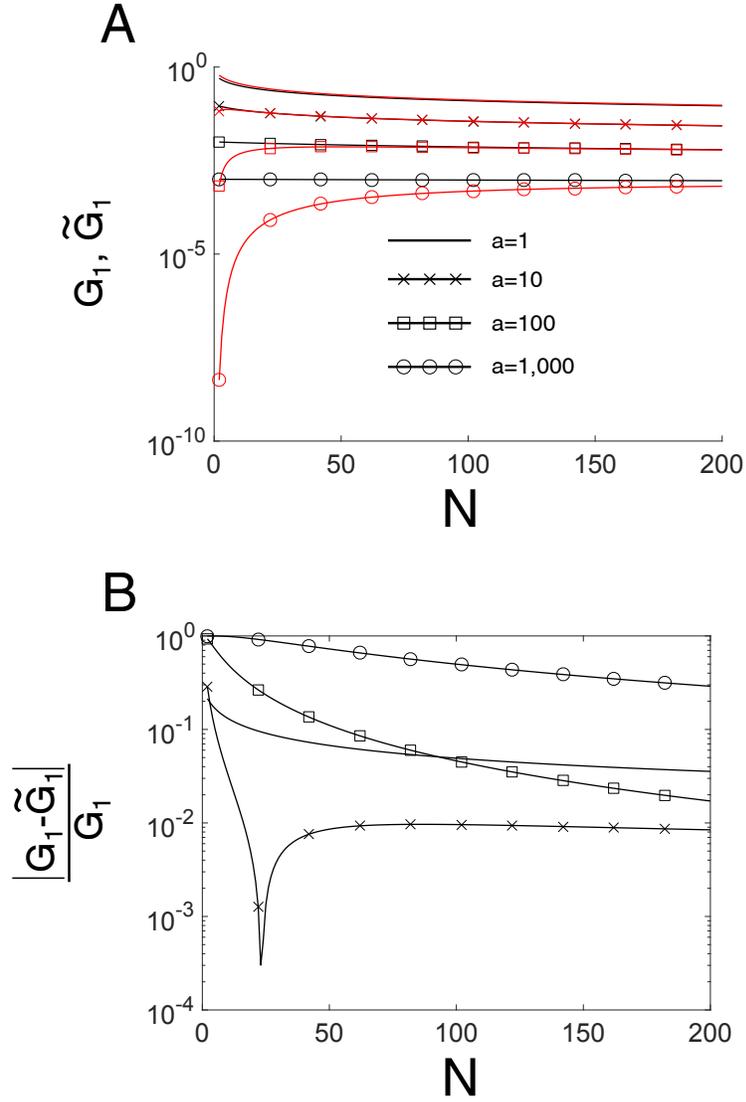}
\caption{{\textbf{Approximation of $G_1$. } (A) Plot of $G_1$ (black) and approximation $\tilde{G}_1$ (red, eq.  \eqref{asymG1})   vs $N$ for  $a=1, 10, 100$ and $1,000$.  (B) Comparison of the  values of  $G_1$ and $\tilde{G}_1$,    measured as $\frac{|G_1-\tilde{G}_1|}{G_1}$.}}
\label{fig:error}
\end{center}
\end{figure}

We provide here approximations for the functions $G_n$  defined in eq. \eqref{definitionGi}. Kummer's function can be expressed in terms of generalized Laguerre polynomials
\beq
_1F_1(-n;b;z)= \frac{n!}{(b)_n}L_{n}^{b-1}(z)
\eeq
where $n$ is an integer. The asymptotic behavior of the Laguerre polynomial for large $n$, fixed $x>0$ and $\alpha$, is given by \cite{Szego1975}
\beq
L_{n}^{\alpha}(-x) \approx  \frac{n^{\frac{\alpha}{2}-\frac1{4}}}{2\sqrt{\pi}} \frac{e^{-\frac{x}{2}}}{x^{\frac{\alpha}{2} +\frac{1}{4}}} \exp \bigg(2\sqrt{x(n + \frac{\alpha+1}{2}) } \bigg).
\eeq
We can now evaluate $G_1$ by finding asymptotic expressions for $_1F_1(-N+1;2;-2a)$ and
 $_1F_1(-N+2;3;-2a)$.
 We have for large $N$,
{ \beq
 _1F_1(-N+1;2;-2a) \approx \frac{1}{N} \frac{(N-1)^{1/4}}{2\sqrt{\pi}} \frac{e^{-a}}{(2a)^{3/4}} \exp \bigg(2\sqrt{2aN } \bigg),  \nonumber \\
  \eeq
 and
  \beq
 _1F_1(-N+2;3;-2a) \approx \frac{1}{N(N-1)} \frac{(N-2)^{3/4}}{2\sqrt{\pi}} \frac{e^{-a}}{(2a)^{5/4}} \exp \bigg(2\sqrt{2a(N-1/2) } \bigg).  \nonumber \\
  \eeq}
Finally
\beq
G_1(a,N)&\approx& \sqrt{\frac{2a}{N}} \exp  \bigg( 2\sqrt{2a(N-1/2)}-2\sqrt{2aN}    \bigg), \nonumber \\
&\approx&   \sqrt{\frac{2}{aN}}\exp  \bigg(- \sqrt{\frac{ a}{2N}}  \bigg) = \tilde{G}_1(a,N) .
\label{asymG1}
\eeq
Similarly, the present computation can be generalized to obtain the asymptotic approximation $\tilde{G}_n$ for the functions $G_n$, valid for large $N$ and fixed $n$,
  \beq
\tilde{G}_n(a,N)&\approx&   \frac{(n+1)!}{(2aN)^{n/2}}\exp  \bigg(- n \sqrt{\frac{ a}{2N}}  \bigg).
\label{eq:MathasymGn}
\eeq
{In Figure \ref{fig:error}, we compare the exact value of $G_1$ (defined in \ref{definitionG1}) with the approximation $\tilde{G}_1$ given by expression \eqref{asymG1}. The approximation is more accurate for intermediates values of $a$ (see Fig. \ref{fig:error}B). In addition, the error function $\frac{|G_1-\tilde{G}_1|}{G_1}$ has a discontinuity in the derivative for $a=10$. Indeed at $a=10$, the function $G_1$ and $\tilde{G}_1$ cross each other and thus $|G_1-\tilde{G}_1|$ has a singular derivative, shown by a cusp type behavior in Fig. \ref{fig:error}B. Moreover, the function $G_1$ is always decreasing with $N$, however its approximation $\tilde{G}_1$ is non monotonic (Fig. \ref{fig:error}A).}
{ Finally, by using the asymptotic expression  \eqref{asymG1}, we find the approximation of the number of clusters for large $N$,
 \beq
\mu_1(a)\approx 1 + \sqrt{2aN}\exp  \bigg(- \sqrt{\frac{ a}{2N}}  \bigg) .
\eeq}
%
%

 {To obtain a numerical approximation of $G_1$ for small $a$,  we used the finite continued fraction decomposition for the ratio of hypergeometric functions $_1F_1$ \cite{Jones}, and obtain
 \beq
G_1(a,N) &=&  \frac{_1F_1(-N+2;3;-2a)}{_1F_1(-N+1;2;-2a)} \nonumber \\
&=&  \cfrac{1}{1+\cfrac{ (N+1)a/3}{1+\cfrac{ (N-2)a/6}{1+\cfrac{ (N+2) a/10}{1+\cfrac{(N-3)a/15}{\ddots+ \frac{ a/(N-1)(2N-3)}{1+a/(N-1)}     }}}}}.
\label{continuedfraction}\eeq
 We obtain the Taylor expansion of $G_1$ for small $a$ from the continued fraction,
 \beq
 G_1(a,N)  &=&1 - \frac{N+1}{3} a+ \bigg(  \frac{(N+1)(N-2) }{18} + \frac{(N+1)^2}{9}\bigg) a^2 \nonumber \\
 & -&  \bigg(   \frac{(N+1)(N-2)(N+2)}{180} + \frac{(N+1)(N-2)^2}{108} + \frac{(N+1)^3}{27}   \bigg)a^3 \nonumber \\
 &+& o(a^3). \nonumber \\
 \eeq
This expression is computed from the continued fraction $G_1$ written
\beq
G_1 (a) & =&  \frac{1}{1+\ds{F_1(a)a}}  \nonumber\\
&=& \frac{1}{1+\ds{\frac{f_1 a}{1+F_2(a) }}},
\eeq
where the general term of the sequence $F_n$ is
\beqq
F_n(a)  =\frac{ f_n }{1+ F_{n+1}(a) a}
\eeqq
with $F_{n+1}(a) =O(1)$. $G_1$ can also be written as
\beq
G_1 (a)  =  \frac{1}{1+\ds{\frac{a f_1}{1+  \ds{ \frac{a f_2 }{1+ \frac{\dots}{1+ a f_n}}} }}}.
\eeq
To obtain a third-order Taylor expansion of $G_1$, we simply used $F_1$, $F_2,$ $F_3$ and we have truncated the rest of continued fraction to obtain
   \beq
 G_1 (a)  &=& 1-F_1(a) a+F_1^2(a)a^2 - F_1^3(a) a^3+o(a^3). \nonumber \\
 \eeq
The expansion of $F_1$ is
\beq
F_1(a)  &=&   \frac{ f_1}{1+F_2(a)a} \nonumber   \\
&=& f_1(  1- F_2(a) a +F_2^2(a)a^2  ) +o(a^2), \nonumber\\
&=& f_1\bigg(  1- f_2(1-f_3 a) a +f_2^2(a)a^2  \bigg) +o(a^2), \nonumber \\
&=& f_1  -f_1f_2 a +(f_1f_2f_3 +f_1f_2^2) a^2 +o(a^2).
\eeq
Using $F_1$ expansion into $G_1$, we finally get
\beq
 G_1 (a)  &=& 1-f_1a  + (f_1^2+f_1f_2)a^2 - (f_1f_2f_3 +f_1f_2^2+f_1^3) a^3    +o(a^3). \nonumber\\
\eeq
}
%
%
%
%
\subsection{Statistics for the number of clusters of a given size}

\begin{figure}[http!]
\begin{center}
\includegraphics[scale=0.65]{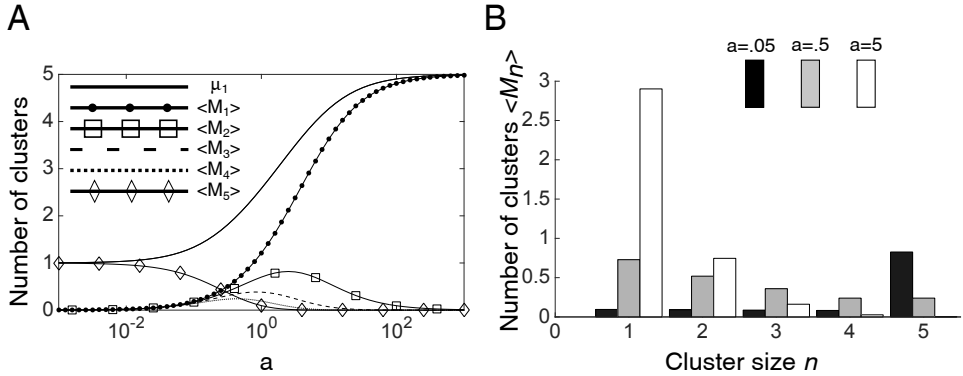}
\caption{ (\textbf{A}) Mean number of clusters of size $n$ as a function of the parameter $a$ for $N=5$ particles (equation \eqref{eq:Mathmeansize}), and total number of clusters $\mu_1$. (\textbf{B}) Mean number of clusters $\langle M_n \rangle$  as a function of the cluster size $n$, for $a=0.05$, $a=0.5$ and $a=5$.}
\label{G1}
\end{center}
\end{figure}

We now compute several moments for the size of clusters when there are $N$ particles. Using relation \eqref{eq:MiNK},  we first obtain the expression for the mean number of clusters of size $n$ when there are $K $ clusters
\beq
\langle M_n \rangle_{N,K}&=&\sum_{(m_i)\in \ANK'}m_np'(m_i|K) \nonumber \\
&=&   a \frac{C_{N-n,K-1}}{C_{N,K}}  . \label{MnCnk}
\eeq
Determining this relation requires   computing the normalizing constant $C_{N,K}$ given in eq. \eqref{cnk}.
The normalizing constant  $C_{N,K}$ is the  $N$th order coefficient of $S^K$,  where $S$ is  the generating function \eqref{generatingfunction}
\beq
S(x)  =  \sum_{i=1}^\infty a_ix^i = a \frac{x}{1-x}.
\eeq
The coefficient $C_{N,K}$  is thus equal to the $N-K$th order coefficient of $\frac{1}{K!}\frac{a^K}{(1-x)^{K}}$ (Remark \ref{remarkgeneratingfunction}). We obtain this coefficient by differentiating $N-K$ times  $\frac{1}{(1-x)^{K}}$ and estimating the derivative at $x=0$. We obtain that
\beq
C_{N,K}& =\frac{1}{K!} & a^K \frac{1}{(N-K)!}K(K+1)\dots(K+N-K-1) \nonumber\\
&=&\frac{a^K}{K!}    \frac{(N-1)!}{(K-1)! (N-K)!}.
\label{CNKaconstant}
\eeq
Thus, by combining \eqref{MnCnk} and \eqref{CNKaconstant},  { we obtain that the number of cluster of size $n$, when the $N$ particles are distributed in $K$ clusters, is}
\beq \label{mjnk}
\langle M_n \rangle_{N,K} = \frac{(N-n-1)!K!(N-K)!}{(N-1)!(K-2)!(N-n-K+1)!},
 \eeq
 which we remark to be independent of $a$.
%
The mean number of clusters of size $n$ is obtained by summing over all {possible} configuration with $K$ clusters,
\beq
 \langle M_n \rangle&=&\sum_{K=1}^N \langle M_n \rangle_{N,K}\Pi_K\\
 &=& \frac{(N-n-1)!}{(N-1)!} \sum_K \frac{K(K-1)(N-K)!}{(N-n-K+1)!} \Pi_K. \nonumber
 \eeq
Using expression \eqref{eq:MathPik} for $\Pi_K$, we obtain
{\beq
\langle M_n \rangle&=&2a\ds{\frac{{_1F_1}(-N+1+n;2;-2a)}{{_1F_1}(-N+1;2;-2a)}} \mbox{ if } n<N,\\
\mbox{and }&& \nonumber \\
 \langle M_N \rangle&=&\ds{\frac{1}{{_1F_1}(-N+1;2;-2a)}}.
\label{eq:Mathmeansize}
\eeq}
The mean number of clusters of size $N$ is exactly equal to the probability $\Pi_{1}(N)$ of having one cluster when there is $N$ particles   (see eq.  \eqref{probability1cluster}). Indeed this is the only configuration where a cluster of size $N$ can appear.
The mean number of clusters of size $n$ is $2a\frac{\Pi_{1}(N)}{\Pi_{1}(N-n)}$, which means that it is given by  the ratio of the probability of having one cluster when there are $N$ particles over the probability of having one cluster when there are $N-n$ particles.\\
The number of clusters can also be written using the function $G_k$ defined in \eqref{definitionGi},
\beq
\langle M_n \rangle&=&\sum_{k=0}^{n}(-1)^k\frac{(2a)^{k+1}}{(k+1)!}\frac{n!}{(n-k)!k!}G_k \\
&=& 2a\sum_{k=0}^{n}(-1)^k\frac{\Pi_{k+1}(n)}{\Pi_1(n)}G_k,
\eeq
where $\Pi_{k}(n)$ is the probability of having $k$ clusters in a system of $n$ particles.
%
To summarize this analysis, we plotted in Fig. \ref{G1} the mean number of clusters of size $n$ for $N=5$ particles.\\

\subsection{Probability to find two particles in the same cluster}
We now evaluate the probability to find two particles in the same cluster for a constant kernel. When there are $N$ particles, we first prove the following
\begin{lemma}  \label{lemmaproba}
\beq
 \sum_{j=1}^{N}j^2 \langle M_j \rangle_{N,K}=N+2N\frac{N-K}{K+1},
\label{sommen2}
\eeq
where
\beq
\langle M_j \rangle_{N,K} = \frac{(N-j-1)!K!(N-K)!}{(N-1)!(K-2)!(N-j-K+1)!}
\eeq
is given by relation \eqref{mjnk}.
\end{lemma}
\Proof
Using formula \eqref{mjnk}, we have
\beq
\sum_{j=1}^{N}j^2 \langle M_j \rangle_N&=& \sum_{j=1}^{N}j^2\frac{(N-j-1)!K!(N-K)!}{(N-1)!(K-2)!(N-j-K+1)!} \nonumber \\
&=& \frac{K(K-1)(N-K)!}{(N-1)!} \sum_{j=1}^{N-K+1}j^2\frac{(N-j-1)!}{(N-j-K+1)!}. \label{eqSNK}
\eeq
To determine the sum in eq. \eqref{eqSNK}, we introduce the sum
\beq
S_{N,K} =\sum_{j=1}^{N-K+1}j^2\frac{(N-j-1)!}{(N-j-K+1)!},
\eeq
and prove now that
\beq
S_{N,K}= \frac{N!}{(N-K)!}\frac{2N-K+1}{(K-1)K(K+1)}.
\eeq
We first obtain a recurrence relation between $S_{N+1,K}$ and $S_{N,K}$
\beq
S_{N+1,K}&=& \sum_{j=1}^{N-K+2}j^2(N-j)(N-1-j)...(N-K+3+j) \nonumber\\
&=&   \sum_{j=0}^{N-K+1}(j+1)^2(N-j-1)(N-j-2)...(N-K+2+j)\nonumber\\
&=& \frac{(N-1)!}{(N-K+1)!} + \sum_{j=1}^{N-K+1}(j^2+2j+1)(N-j-1)(N-j-2)...(N-K+2+j) \nonumber \\
&=& \frac{(N-1)!}{(N-K+1)!}  + S_{N,K}  + A+2B,  \nonumber \\ &&
\eeq
where
\beq
A=   \sum_{j=1}^{N-K+1}(N-j-1)(N-j-2)...(N-K+2+j)
\eeq
and
\beq
B=  \sum_{j=1}^{N-K+1}j(N-j-1)(N-j-2)...(N-K+2+j).
\eeq
We first use a change of variable for $A$ and find
\beq
A&=&\sum_{j=K-2}^{N-2} j(j-1)...(j-K+3) \nonumber\\
&=&\sum_{j=K-2}^{N-2} \frac{j!}{(j-K+2)!} \nonumber\\
&=& (K-2)!\sum_{j=K-2}^{N-2} \binom{j}{K-2} \nonumber.
\eeq
Using the binomial relation
\beq
\sum_{j=k}^{n}\binom{j}{k} = \binom{n+1}{k+1} \nonumber,
\eeq
we obtain
\beq
A&=& (K-2)!\binom{N-1}{K-1} \nonumber\\
&=&\frac{(N-1)!}{(N-K)!(K-1)}.
\eeq
Similarly, we obtain for $B$
\beq
B&=&\sum_{j=K-2}^{N-2}(N-j-1)\frac{j!}{(j-K+2)!} \nonumber\\
&=&NA- \sum_{j=K-2}^{N-2}\frac{(j+1)!}{(j-K+2)!} \nonumber\\
&=& NA - (K-1)!\sum_{j=K-1}^{N-1} \binom{j}{K-1} \nonumber\\
&=& \frac{N!}{(N-K)!(K-1)}  - \frac{N!}{(N-K)!K} \nonumber\\
&=&\frac{N!}{(N-K)!(K-1)K}.
\eeq
Finally we obtain the induction relation for $S_{N,K}$
\beq
S_{N+1,K}&=&S_{N,K} + 2\frac{N!}{(N-K)!(K-1)K}+\frac{(N-1)!}{(N-K)!(K-1)} + \frac{(N-1)!}{(N-K+1)!} \nonumber\\
&=& S_{N,K} + \frac{N!(2N-K+2)}{(N-K+1)!(K-1)K}.
\eeq
Using that $S_{K,K} = (K-2)! $, we finally evaluate the sum
\beq
S_{N,K }&=& \frac{1}{(K-1)K}\sum_{j=K}^{N} \frac{(j-1)!(2j-K)}{(j-K)!} \nonumber\\
&=&  \frac{2}{(K-1)K}\sum_{j=K}^{N}\frac{j!}{(j-K)!} -\frac{1}{(K-1)!}\sum_{j=K}^{N}\frac{(j-1)!}{(j-K)!} \nonumber\\
&=& 2(K-2)!\sum_{j=K}^{N}\binom{j}{K}-(K-2)!\sum_{j=K-1}^{N-1}\binom{j}{K-1} \nonumber\\
&=& 2(K-2)!\binom{N+1}{K+1}-(K-2)!\binom{N}{K} \nonumber\\
&=& \frac{(K-2)!}{(N-K)!}\big(2 \frac{(N+1)!}{(K+1)!}- \frac{N!}{K!}\big) \nonumber\\
&=&\frac{N!}{(N-K)!}\frac{2N-K+1}{(K-1)K(K+1)},
\eeq
which is formula \eqref{sommen2}. \Endproof

\begin{thm}
The probability to find two particles in the same cluster is
\beq
\langle P_2 \rangle  =G_1,
\eeq
where $G_1$ is defined in \eqref{definitionG1}.
\end{thm}
\Proof
Using   eq. \eqref{eq:P2generalcase} and lemma \ref{lemmaproba}, we can now compute the probability to find two particles in the same cluster
\beq
\langle P_2 \rangle &=& \frac{1}{N(N-1)}\sum_{K=1}^N\Pi_K
 \sum_{j=1}^{N}j^2 \langle M_j \rangle_{N,K}-\frac{1}{N-1} \nonumber \\
 &=& \frac{1}{N(N-1)}\sum_{K=1}^N\Pi_K
\left( N+2N\frac{N-K}{K+1} \right)-\frac{1}{N-1}.
\eeq
Thus,
\beq
\langle P_2 \rangle&=&\frac{2}{N-1}\sum_{K=1}^N\Pi_K\frac{N-K}{K+1}  \nonumber\\
&=&-\frac{2}{N-1}+2\frac{N+1}{N-1}\sum_{K=1}^N\frac{1}{K+1}\Pi_K.
\label{eqP2_1}
\eeq
Following formula \eqref{deriv1}, the sum in eq. \eqref{eqP2_1} can be expressed by integrating the Kummer's function
{\beq
\sum_{K=1}^N\frac{1}{K+1}\Pi_K = \frac{\Pi_1}{2}\left. \frac{\int z _1F_1(-N+1;2;z)dz}{z^2} \right|_{z=2a}.
\eeq}
Integrating the hypergeometric series gives
\beq
\int z {_1F_1}(-N+1;2;z) dz=z^2  {_2F_2}(-N+1,2;2,3;z),
\eeq
which leads to
\beq
\sum_{K=1}^N\frac{1}{K+1}\Pi_K &=& \frac{1}{2}\frac{_2F_2(-N+1,2;2,3;-2a)}{_1F_1(-N+1;2;-2a)} \nonumber \\
&=& \frac{1}{2}\frac{_1F_1(-N+1;3;-2a)}{_1F_1(-N+1;2;-2a)}.
\label{eqP2_2}
\eeq
By combining eq. \eqref{eqP2_1} and  \eqref{eqP2_2}
we obtain
\beq
\langle P_2 \rangle=-\frac{2}{N-1} + \frac{N+1}{N-1}\frac{_1F_1(-N+1;3;-2a)}{_1F_1(-N+1;2;-2a)}.
\eeq
The three-term recurrence relation for Kummer's function (\cite{Abramowitz} Eqs. 13.4.1Ð-13.4.6) gives
\beq
_1F_1(-N+1;3;-2a)=   \frac{N-1}{N+1}{_1F_1}(-N+2;3;-2a)   + \frac{2}{N+1}{_1F_1}(-N+1;2;-2a).  \nonumber \\
\eeq
Finally, using eq. \eqref{definitionG1}, we obtain
\beq
\langle P_2 \rangle=  G_1. \mbox{ \Endproof }
\label{eq:MathP2model2}
\eeq
To finish this section, we note that the large $N$ asymptotic of the probability that two particles are in the same cluster is
\beq
\langle P_2 \rangle\approx  \sqrt{\frac{2}{aN}}.
\eeq
{ Many results presented in this section can be used to study the distribution of clusters in biological systems such as telomere organization in yeast. We provided here the explicit derivations of the exact and asymptotic formulas that can be used to analyze experimental and simulation results \cite{Hoze2012}.}

\section{Example 2: the case $a_i=a$ for $i<M$ and $a_i=0$ if $i\geq M$} \label{section:example2}
In this section, we consider $N$ particles that can associate or dissociate at a constant rate, but in addition they cannot form clusters of more than $M$ particles. The configuration space for distributions of $N$ particles in $K$ clusters of size less than $M$ is
\beq
\ANKM'= \{ (m_i)_{1\leq i \leq M} ; \sum_{i=1}^M im_i=N \mbox{, } \sum_{i=1}^{M}
m_i=K\}.
\eeq
First, the minimal number of clusters is necessarily bounded by $K\geq N/M$, since the opposite would imply a cluster of at least $M+1$ particles. The probability of a configuration $(m_1,...,m_M) \in \ANKM'$ is equal to
\beq
Pr\{(m_1,...,m_M) \in \ANKM'\} =\frac{1}{C_{N,K,M}}\frac{1}{m_1!...m_M!},
\eeq
where the normalization constant $C_{N,K,M}$ is the $N$th order coefficient of
\beq
(aX+aX^2+...+aX^M)^K&=&a^KX^K(\frac{X^M-1}{X-1})^K \nonumber \\
&=&a^K(\frac{X}{1-X})^K\sum_{n=0}^{K}\binom{K}{n}(-1)^nX^{nM} \\
&=&a^K \frac{1}{(1-X)^K}\sum_{n=0}^{K}\binom{K}{n}(-1)^nX^{nM+K}.
\eeq
Then the $N$th order coefficient of the polynomial is obtained by finding the $(N-nM-K)$th order coefficient of $(1-X)^{-K}$
\beq
C_{N,K,M}   = a^K  \sum_{n=0}^{K}\binom{K}{n}(-1)^n\frac{1}{(N-(nM+K))!}D^{(N-(nM+K))}\left( \frac{1}{(1-X)^K}\right)_{|X=0}, \nonumber \\
 \eeq
where we write $D^{(n)}$ the $n-$th order derivative. Thus, setting $K_0 =\lfloor \frac{N-K}{M} \rfloor$, where $\lfloor . \rfloor$ is the floor function,  we have
\beq
C_{N,K,M}  = a^K K\sum_{n=0}^{K_0} \frac{(N-nM-1)!}{n!(K-n)!(N-(nM+K))!}(-1)^n    .
\label{limitsizeCNKM}
\eeq
For $M=N$ we find $K_0=0$ and the normalization constant
\beq
C_{N,K,N}= a^K\frac{(N-1)!}{(K-1)!(N-K)!},
\eeq
is equal to the normalization constant $C_{N,K}$ obtained for the constant kernel in section \ref{section:example1}.\\
The mean number of clusters of size $i\leq M$ conditioned on the number of clusters $K$ is
\beq
\langle M_i \rangle_K  &=& \sum_{m_i\in \ANKM'} m_ip(m_1, ...m_M) \nonumber \\
&=& a\frac{C_{N-i,K-1,M}}{C_{N,K,M}}.
\label{limitsizeMi}
\eeq
To find the probability to have $K$ clusters, we now redefine the formation rate. In section \ref{section:example1}, the formation rate was proportional to the number of pairs of particles since all of them could form a new cluster. In the present case, two clusters of size $i$ and $j$ can form a new cluster only if $i+j\leq M$. The formation rate when there are $K$ clusters is thus
\beq
f_K = \sum_{(m_i)\in \ANKM'}  p\mn \left( \sum_{i=1}^{M/2}\frac{m_i(m_i-1)}{2}+\sum_{\underset{ i+j\leq M; i \neq j}{i,j=1}}^{M}m_im_j \right).  \nonumber \\
\eeq
The formation rate can be written as a function of the coefficients $C_{N,K,M}$ as
\beq
f_2 = C_{N,2,M},
\eeq
and for $K>2$
\beq
f_K &=&   \frac{K(K-1)}{2}\sum_{i=1}^{\min( \frac{M}{2} , \frac{N-K+2}{2})}C_{N-2i,K-2,M} \nonumber\\
&+&\frac{K(K-1)}{2}\sum_{\underset{ i+j\leq M}{i,j=1}}^{\min(M-1,N-K+1)}C_{N-i-j,K-2,M}.
\eeq
The separation rate remains unchanged $s_K=(N-K)a$, and the probabilities at steady state are given by
\beq
\Pi_K = \frac{f_{K+1}}{s_K}\Pi_{K+1}.
\eeq
A simple expression is certainly hopeless, but the limit  $a\rightarrow0$ is informative be: contrary to the previous case with a constant kernel ($M=N$), where the particles form a single cluster of size $N$, the present system contains multiple steady state distributions. The clusters grow independently and reach either their limit size $M$ or are configured such that the sum of the sizes of each pair of cluster is larger than $M$. All possible configurations contain exactly $\lceil N/M\rceil$ clusters, where $\lceil . \rceil$ is the ceiling function.

We illustrate the limit case $a\rightarrow 0 $ for $N=9$, $M=4$ (Fig. \ref{coaglimitsize}). Because $a>0$, all partitions  are accessible, but as $a\rightarrow 0$, the steady state configurations are dominated by the configurations with the largest possible cluster size $(4,4,1)$, $(4,3,2)$ and $(3,3,3)$.
Applying formulas \eqref{limitsizeCNKM} and \eqref{limitsizeMi}, we obtain the limit cluster configuration  probabilities
\beq
p(4,4,1) &=& \frac{3}{10} \nonumber \\
p(4,3,2) &=& \frac{6}{10} \label{configurationM4} \\
p(3,3,3) &=& \frac{1}{10} . \nonumber
\eeq
{These steady state probabilities do not depend on the initial particles configurations as long as $a\neq 0$. { For $a=0$, there are three possible configurations $(4,4,1)$, $(4,3,2)$ and $(3,3,3)$: once equilibrium is attained, the clusters will remain unchanged. The probability to get to equilibrium depends on the configuration and the order of clustering events. } When there is no limitation in the cluster formation ($M=N=9$), a single cluster containing all particles is formed (Fig. \ref{coaglimitsize}, left panel). For large values of $a$, most clusters are very small, and the distributions are similar for $M=4$ and $M=9$ (Fig. \ref{coaglimitsize}, right panel).
\begin{figure}[http!]
\begin{center}
\includegraphics[scale=0.6]{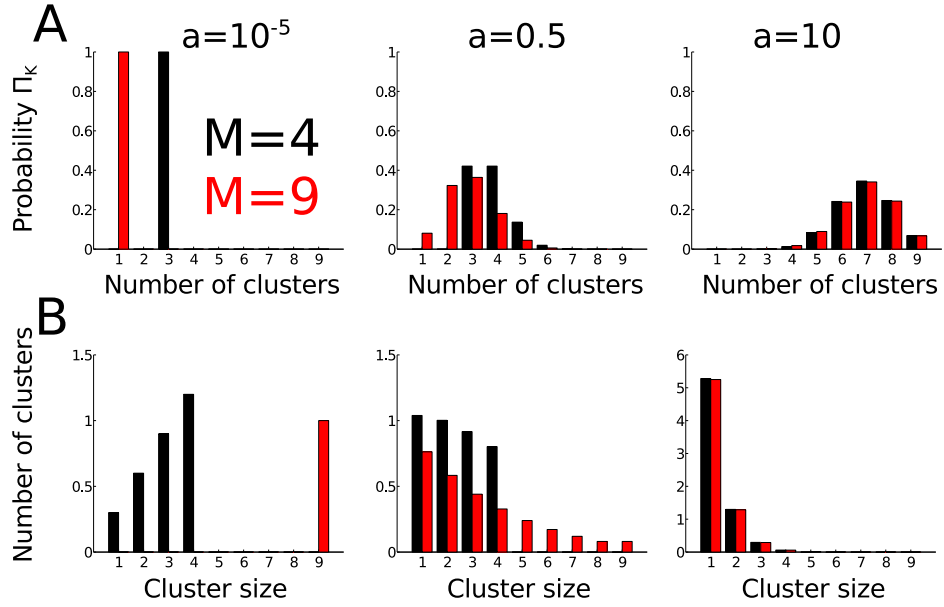}
\caption{ (\textbf{A}) Distribution of the number of clusters $\Pi_K$ for $N=9$, when cluster sizes are limited ($M=4$, black) and not limited  ($M=9$, red).  There is a minimum of $\lceil N/M\rceil$ clusters. From left to right : $a=10^{-5}$, $a=0.5$, $a=10$.
(\textbf{B})   Mean number of clusters of each size $\langle M_n \rangle$. For $a\rightarrow 0$, for $N=9$ and $M=4$ the clusters  organize in three different cluster configurations, while for $M=N$ a single cluster containing $N$ particles is formed.  }
\label{coaglimitsize}
\end{center}
\end{figure}

The probability for two particles to be in the same cluster provides a good estimation for the cluster distribution for various values of the parameter $a$ (Fig.  \ref{fig:p2}). When $a$ is large, most particles are contained in very small clusters and the probability $\langle P_2 \rangle$ is similar for the cases $M=4$ and $M=9$. When $a\rightarrow 0$, particles tend to form larger clusters. A single cluster containing all particles is formed and $\langle P_2 \rangle\rightarrow 1$ when $M=9$,  but the maximal value of $\langle P_2 \rangle$ is less than 1 when the maximal cluster size is limited. We can explicitly compute $\langle P_2 \rangle$ in the limit case $a\rightarrow 0$. For example for $M=4$, using eq. \eqref{Matheqq}, and summing over all possible configurations \eqref{configurationM4}, we obtain
\beqq
\langle P_2 \rangle &=& p(4,4,1) P_2(4,4,1) + p(4,3,2)P_2(4,3,2) +p(3,3,3) P_2(3,3,3) \\
&=&  \frac{3}{10} \frac{24}{72} + \frac{6}{10} \frac{ 20}{72} + \frac{1}{10} \frac{18}{72} \\
&=& \frac{7}{24}.
\eeqq



\begin{figure}[http!]
\begin{center}
\includegraphics[scale=0.4]{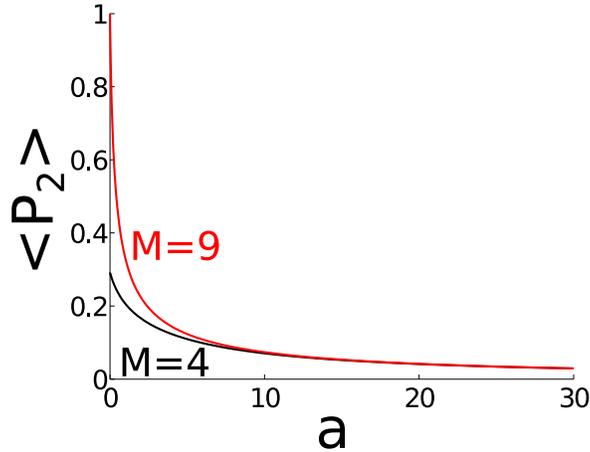}
\caption{ {\bf Probability  $\langle P_2 \rangle$  that two particles are in the same cluster.} The parameters are  $N=9$ and $M=4$ (black), $M=9$ (red). For large values of $a\gg1$, only small clusters are present and the steady state distributions are similar for the cases $M=4$ and $M=9$. When $a\rightarrow0$ the clusters organize in three different cluster configurations, while for $M=N$ a single cluster containing $N$ particles is formed.}
\label{fig:p2}
\end{center}
\end{figure}

\section{Example 3: Application to the case $a_i=ai$} \label{section:example3}
We finally consider the case $a_i=ai$.
It has been shown \cite{Durrett1999} that the number of clusters of size $i$ is asymptotically
\beq
\langle M_i \rangle =  aie^{-i\sqrt{2a/N}}.
\eeq
The generating function $S$ (eq. \eqref{generatingfunction}) is given by
\beq
S(x)   =a \frac{x}{(1-x)^2},
\eeq
which gives that
\beq
C_{N,K}  & =& a^K\frac{1}{(N-K)! } D^{N-K}\frac{1}{(1-x)^{2K}}_{|x=0} \nonumber \\
&=&a^K\frac{1}{(N-K)! } \frac{(N+K-1)!}{(2K-1)! } \nonumber \\
&=&a^K \binom{N+K-1}{N-K}.
\eeq
We thus obtain using formula \eqref{Migeneralcase} that
\beq
\langle M_i \rangle_{N,K} =  i \frac{ \binom{N-i+K-2}{N-K-i+1}}{ \binom{N+K-1}{N-K}  }.
\eeq
To obtain the number of clusters of size $i$, we determine the probability of a distribution of $K$ clusters $\Pi_K$. We consider the   coagulation kernel $C(i,j)=1$ and the fragmentation kernel $F(i,j) = a \frac{ij}{i+j} $, and obtain that
\beq
d(n) &=& \sum_{i=1}^{n-1}a\frac{i(n-i)}{n} \nonumber\\
&=& \frac{a(n^2-1)}{6}.
\eeq
The separation rates are
\beq
s_1=\frac{a(N^2-1)}{6}
\eeq
and for $K\geq 2$
\beq
s_K&=& \frac{\sum_{i=1}^{N-K+1}d(i)a_i C_{N-i,K-1}}{C_{N,K}} \nonumber\\
&=& \frac{a}{6} \frac{\sum_{i=1}^{N-K+1}i(i^2-1) \binom{N-i+K-2}{N-i-K+1}}{\binom{N+K-1}{N-K}}\nonumber \\
&=& \frac{a}{6}\frac{1}{\binom{N+K-1}{N-K}}\frac{1}{(2K-3)!} \sum_{i=1}^{N-K+1} \frac{i(i^2-1)(N-i+K-2)!}{(N-i-K+1)!}.\label{expressionfaibi}
\eeq
To go further in the determination of $s_K$, we evaluate the sum
\beq
\varphi_{N,K}=\frac{1}{(2K-3)!}\sum_{i=1}^{N-K+1} \frac{(i^3-i)(N-i+K-2)!}{(N-i-K+1)!}.
\eeq
After the change of variables $j=N-i-K+1$ in the sum we obtain
\beq
\varphi_{N,K}&=&\frac{1}{(2K-3)!}\sum_{j=0}^{N-K}\left((N-K-j+1)^3-(N-K-j+1)\right)\frac{(j+2K-3)!}{j!}\nonumber\\
&=& \sum_{j=0}^{N-K}\left((N-K-j+1)^3-(N-K-j+1)\right)\binom{j+2K-3}{2K-3}\nonumber\\
&=& \sum_{j=2K-3}^{N+K-3}\left((N+K-2-j)^3-(N+K-2-j)\right)\binom{j}{2K-3}. \nonumber\\
\eeq
We expand the sum and write $\varphi_{N,K}$, with $G(N,K)=N+K-2$, as
\beq
\varphi_{N,K}&=& (G(N,K)^3-G(N,K)) \sum_{j=2K-3}^{G(N,K)-1}\binom{j}{2K-3} \nonumber\\
&+& (1-3G(N,K))\sum_{j=2K-3}^{G(N,K)-1}j\binom{j}{2K-3}\nonumber \\
&+&3G(N,K) \sum_{j=2K-3}^{G(N,K)-1}j^2\binom{j}{2K-3}\nonumber\\
&-&  \sum_{j=2K-3}^{G(N,K)-1}j^3\binom{j}{2K-3}.
\eeq
We evaluate the sums using the formulas
\beq
\sum_{j=k}^{n}\binom{j}{k}&=&\binom{n+1}{k+1}\\
\sum_{j=k}^{n}j\binom{j}{k}&=&(k+1)\binom{n+2}{k+2}-\binom{n+1}{k+1}\\
\sum_{j=k}^{n}j^2\binom{j}{k}&=&(k+1)(k+2)\binom{n+3}{k+3}   \\ &-&3(k+1)\binom{n+2}{k+2}+\binom{n+1}{k+1}\nonumber\\
\sum_{j=k}^{n}j^3\binom{j}{k}&=&(k+1)(k+2)(k+3)\binom{n+4}{k+4}    \\ &-&  6(k+1)(k+2)\binom{n+3}{k+3}  \nonumber \\
&+& 7(k+1)\binom{n+2}{k+2}-\binom{n+1}{k+1} \nonumber .
\eeq
So the sum $\varphi_{N,K}$ is equal to
\beqq
\varphi_{N,K}&=&(G(N,K)^3-G(N,K))\binom{G(N,K)}{2K-2}+(1-3G(N,K))(2K-2)\binom{G(N,K)+1}{2K-1}
\\&-&(1-3G(N,K))\binom{G(N,K)}{2K-2}+3G(N,K)(2K-2)(2K-1)\binom{G(N,K)+2}{2K}\\
&-&9(2K-2)G(N,K)\binom{G(N,K)+1}{2K-1}+3G(N,K)\binom{G(N,K)}{2K-2}\\
&-&(2K-2)(2K-1)2K\binom{G(N,K)+3}{2K+1}+6(2K-2)(2K-1)\binom{G(N,K)+2}{2K}\\
&-&7(2K-2)\binom{G(N,K)+1}{2K-1}+\binom{G(N,K)}{2K-2},
\eeqq
which can be simplified into
\beq
\varphi_{N,K}&=&\binom{G(N,K)}{2K-2} (G(N,K)^3+5G(N,K))\nonumber\\
&-&6\binom{G(N,K)+1}{2K-1}(  1+2G(N,K))(2K-2)       \nonumber  \\
 &+&3\binom{G(N,K)+2}{2K}(G(N,K)+2)(2K-2)(2K-1)\nonumber\\
&-&\binom{G(N,K)+3}{2K+1}(2K-2)(2K-1)2K.
\eeq
We now use that  $\binom{n+1}{k+1}=\frac{n+1}{k+1}\binom{n}{k}$ to write $\varphi_{N,K}$ as a function of $\binom{G(N,K)+1}{2K-1} $. Using eq. \eqref{expressionfaibi} we obtain
\beq
s_K&=&\frac{a}{6} \frac{2K-1}{G(N,K)+1} (G(N,K)^3+5G(N,K))\nonumber\\
&-&a(  1+2G(N,K))(2K-2)       \nonumber  \\
 &+&\frac{a}{2} (G(N,K)+2)^2\frac{(2K-2)(2K-1)}{2K}\nonumber\\
&-&\frac{a}{6} (G(N,K)+2)(G(N,K)+3)\frac{(2K-1)(2K-2)}{2K+1}.
\eeq
The formation rates are given by
\beq
f_K=\frac{K(K-1)}{2}
\eeq
and the probability of having $K$ clusters is given by the  relation
\beq
\Pi_K = \frac{f_{K+1}}{s_K}\Pi_{K+1}.
\eeq
%
%
%

\section{Conclusion}
In this paper, we investigated a certain class of discrete coagulation-fragmentation processes with a finite number of particles. We determined the steady state probability distribution when the number of clusters is fixed. We studied the cluster distributions using the partitions of the total number of particles with a given number of clusters. We computed the distribution probability function in terms of multinomial coefficients.\\
This approach allows computing various statistical quantities and moments, including the mean number of clusters of a given size conditioned on the total number of clusters. However, computing other quantities, such as the size of the largest clusters cannot be derived from the present results and requires novel methods of calculation. \\
Finally, we defined two new times to characterize the cluster dynamics: one is the time that two particles spend together and second is the time they spend separated. We computed here the fraction of these times, which is the probability that two particles are in the same cluster. We have applied these results to specific coagulation-fragmentation kernels. For the constant kernel, we obtained exact expressions of the number of clusters in terms of hypergeometric function. When the size of the cluster is limited, we obtained a model of nucleation in the limit $a\rightarrow 0$ and found multiple steady state distributions, depending on the initial number of particles and the limit size.

{Our study  on coagulation-fragmentation of a finite number of particles was motivated by stochastic processes  in chemical reaction theory and in molecular and cell biology of the cell nucleus organization.}
When the clusters and free particles evolve in a homogeneous region in dimension 2 or 3, the time two particles spend separated (recurrence time) is exactly the reciprocal of the forward rate of a chemical reaction. However, when the region is not homogenous so that clustering can occur  preferentially in some subregions, this is not anymore the case, and the recurrence time can be shorter than the meeting time, as discussed in the context of telomere clustering in yeast \cite{Hoze2013}. In that case, clustering favors encounter. The time that two particles stay  in the same cluster is an indicator of the possible exchange of genetic information between clustered telomeres, a process that remains to be studied both experimentally and theoretically.

\end{document}